\documentclass[ijoc,sglanonrev]{informs4}

\RequirePackage{tgtermes}
\RequirePackage{newtxtext}
\RequirePackage{newtxmath}
\RequirePackage{bm}
\RequirePackage{endnotes}

\usepackage[preprint, nonatbib]{}

\usepackage{graphicx}
\usepackage{amsmath}
\usepackage{amsfonts}
\usepackage{xcolor}
\usepackage{algorithm}
\usepackage{algpseudocode}
\usepackage{tikz}
\usepackage{caption}
\usepackage{subcaption}
\usepackage{multirow}
\usepackage{makecell}
\usepackage{bm}

\usepackage{natbib}
 \bibpunct[, ]{(}{)}{,}{a}{}{,}%
 \def\bibfont{\small}%

 \EquationsNumberedThrough
 \TheoremsNumberedThrough
\ECRepeatTheorems  %

\MANUSCRIPTNO{IJOC-0001-2024.00}

\newcommand{\CPS}{PCAB~}

\begin{document}

\RUNAUTHOR{Barbato, Ceselli and Messana}
\RUNTITLE{Mathematical programming algorithms for convex hull approximation with a hyperplane budget}

\TITLE{Mathematical Programming Algorithms for Convex Hull Approximation with a Hyperplane Budget}

\ARTICLEAUTHORS{%
\AUTHOR{Michele Barbato, Alberto Ceselli, Rosario Messana}
\AFF{Department of Computer Science, Università degli Studi di Milano, Via Giovanni Celoria 18, 20133 Milan, Italy \\
\textbf{Contact: }\EMAIL{michele.barbato@unimi.it}, \EMAIL{alberto.ceselli@unimi.it}, \EMAIL{rosario.messana@unimi.it}}}

\ABSTRACT{
    We consider the following problem in computational geometry: given, in the $d$-dimensional real space, a set of points marked as positive and a set of points marked as negative, such that the convex hull of the positive set does not intersect the negative set, find $\mathcal{K}$ hyperplanes that separate, if possible, all the positive points from the negative ones. That is, we search for a convex polyhedron with at most $\mathcal{K}$ faces, containing all the positive points and no negative point. The problem is known in the literature for pure convex polyhedral approximation; our interest, stems from its possible applications in constraint learning, where points are feasible or infeasible solutions of a Mixed Integer Program, and the $\mathcal{K}$ hyperplanes are linear constraints to be found.

    We cast the problem as an optimization one, minimizing the number of negative points inside the convex polyhedron, whenever exact separation cannot be achieved. We introduce models inspired by support vector machines and we design two mathematical programming formulations with binary variables. We exploit Dantzig-Wolfe decomposition to obtain extended formulations, and we devise column generation algorithms with ad-hoc pricing routines.
    
    We compare computing time and separation error values obtained by all our approaches on synthetic datasets, with number of points from hundreds up to a few thousands, showing our approaches to perform better than existing ones from the literature. Furthermore, we observe that key computational differences arise, depending on whether the budget $\mathcal{K}$ is sufficient to completely separate the positive points from the negative ones or not. On $8$-dimensional instances (and over), existing convex hull algorithms become computational inapplicable, while our algorithms allow to identify good convex hull approximations in minutes of computation. They can also be used for parametric analyses, being the hyperplane budget a hyperparameter trading approximation quality for model complexity.
}


\maketitle

\section{Introduction}

We consider the following geometric problem. Two finite sets of points in $\mathbb R^d$ are given, labelled as positive and negative respectively. We assume that the convex hull of the positive points does not contain any negative point. Our aim is to find a convex polyhedron having at most $\mathcal{K}$ faces, which contains all the positive points and, if possible, no negative points. \\

Besides its applications in fields like image processing \citep*{jayaram2016convex}, and the interest into the combinatorial problem on its own, we are motivated by its perspective use in constraint learning \citep*{lombardi2018boosting}. Let us assume that a constrained optimization problem is given, together with a suitable sampling of solutions (encoded as points in $\mathbb R^d$), some of them being feasible (the positive ones) and others being infeasible (the negative ones). Let us also assume to know that the region of feasible solutions can be described as a polyhedron, defined by at most $\mathcal{K}$ linear constraints. The problem of finding them is equivalent to our problem described above: finding a convex polyhedron having at most $\mathcal{K}$ faces, enclosing all the points that encode feasible solutions, and excluding all those encoding infeasible ones. Hardness results date back to \citet{megiddo88}.

Without the requirement of limiting the number of faces to $\mathcal{K}$, the problem could be solved by producing the convex hull of the feasible points using well established techniques \citep*{Barber1996TheQA}, and in particular by finding the equations describing the facets of the convex hull. The set of negative points becomes irrelevant, but the combinatorial complexity of these approaches makes them computationally viable only for positive sets of small size. Faster alternatives exist for convex hull approximation \citep*{Balestriero2022DeepHullFC}, adjusting the solution in post processing to contain all the positive points, when necessary.\\

Since we assume that the budget $\mathcal{K}$ is given as input, it is possible that $\mathcal{K}$ hyperplanes are not sufficient to separate all the positive points from the negative ones. Our goal, in this case, is to produce a \emph{small} approximation of the convex hull. We therefore map the problem to an optimization one, in which the number of negative points falling inside the convex polyhedron is minimized (and no positive points falls outside the polyhedron). The intuition behind this choice is the following: if the points in the positive and negative sets came from a uniform random drawing (and were dense enough), then the $\mathcal{K}$ hyperplanes define a convex polyhedron that over-approximates the convex hull of the feasible points, and whose volume is minimized.
\begin{figure}
  \begin{subfigure}{0.31\textwidth}
    \includegraphics[width=\linewidth]{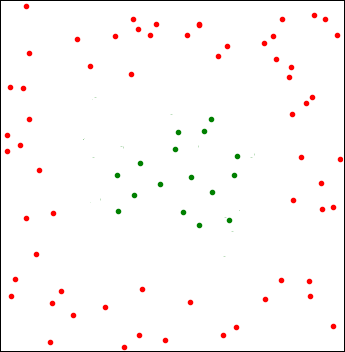}
    \caption{Dataset} \label{fig:1a}
  \end{subfigure}%
  \hspace*{\fill}   
  \begin{subfigure}{0.31\textwidth}
    \includegraphics[width=\linewidth]{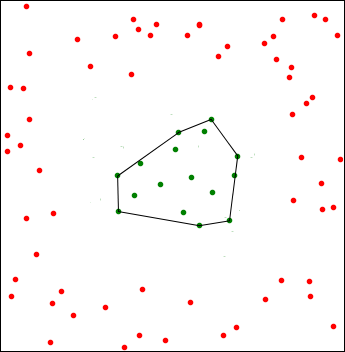}
    \caption{Convex hull} \label{fig:1b}
  \end{subfigure}%
  \hspace*{\fill}   
  \begin{subfigure}{0.31\textwidth}
    \includegraphics[width=\linewidth]{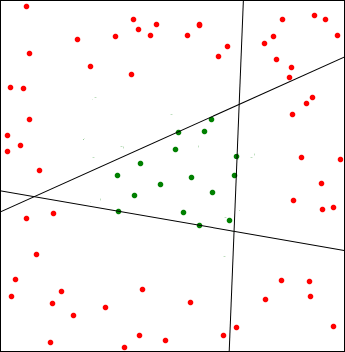}
    \caption{Hull approximation} \label{fig:1c}
  \end{subfigure}
\caption{Given a dataset of positive and negative points (Fig.~\ref{fig:1a}), one can alternatively compute the convex hull of the positive points (Fig.~\ref{fig:1b}) or fix a budget on the number of hyperplanes and compute an over-approximation of the convex hull that satisfies the budget requirement and does not contain any of the negative points (Fig.~\ref{fig:1c}).}
\label{fig:1}
\end{figure}

\paragraph{Contributions.}
We formally introduce the problem and we discuss similarities and differences with related applications from the literature (Section \ref{sec:ProblemStatement}).
Then, we investigate the use of mixed-integer linear programming (MILP) to design new resolution algorithms. In particular, we propose two mixed-integer linear programming models (Section \ref{sec:CompactModels}). We derive extended formulations from them, through the use of Dantzig-Wolfe decomposition.
We design column generation algorithms, employing techniques such as early stopping, and pricing solutions polishing. Additionally, we design a heuristic parallel algorithm for the pricing problem that proves to be empirically advantageous. We also introduce fast upper bounding and master initialization greedy methods (Section \ref{sect:algorithms}).
We compare with both state-of-the-art convex hull algorithms, and to alternative polyhedral separation techniques from the literature on a wide range of instances (Section \ref{sec:ExperimentalSettingAndResults}), focusing on how dimensionality, number of points and hyperplane budget impact the behaviour of each method. Some conclusions are finally drawn (Section \ref{sec:Conclusions}), highlighting how our methods improve on existing ones.

\section{Notation and Background}
\label{sec:ProblemStatement}

We call our problem \emph{polyhedral convex hull approximation with hyperplane budget} (PCAB). It can be formalized as follows. We are given two indexed sets $\mathcal{D}^+ = \{\boldsymbol{a}^+_i\}_{i=1}^m$ and $\mathcal{D}^- = \{\boldsymbol{a}^-_i\}_{i=1}^n$ of points in $\mathbb{R}^d$ 
and a positive integer $\mathcal{K}$, the \emph{hyperplane budget}. \\

A hyperplane $h \in \mathbb{R}^d$ is called \emph{valid} if it admits a valid equation, i.e. it exist $b \in \mathbb{R}$ and $\boldsymbol{w} \in \mathbb{R}^d$ such that $h=\{\boldsymbol x\in\mathbb R^d\colon b+\boldsymbol w\cdot \boldsymbol x=0\}$ and for each point $a_i \in \mathcal{D}^+$ it holds that $b + \sum_{j=1}^d w_j a_{ij} \geq 0$. For the sake of brevity, we identify an equation with its coefficients $(b, \boldsymbol{w})$. \\

A valid solution to \CPS is a set of $\mathcal{K}$ valid equations $s = \left\{(b^k, \boldsymbol{w}^k)\right\}_{k=1}^{\mathcal{K}}$. The \emph{separation error} associated with $s$ is:
\begin{center}
    $err(s) = \left| \left\{ \boldsymbol{a}^-_i\ |\ i \in \{1, \ldots, n\} \wedge b^k + \sum_{j=1}^d w_j^k a_{ij} \geq 0\ \forall k \in \{1, \ldots, \mathcal{K}\} \right\} \right|$.
\end{center}
A valid solution $s^*$ is optimal if $err(s^*)$ is minimum over all valid solutions. \CPS consists in the problem of finding one of these optimal solutions.

\noindent
{\bf Terminology. }
We define the following index sets: $I^+ := \{1, \ldots, m\}$, $I^- := \{1, \ldots, n\}$, $I := \{1, \ldots, m + n\}$, $J := \{1, \ldots, d\}$, $K := \{1, \ldots, \mathcal{K}\}$. \\

Let $\ell: \mathcal{D}^+ \cup \mathcal{D}^- \rightarrow \{-1, 1\}$ be the labeling function such that $\ell(\mathcal{D}^+) = 1$ and $\ell(\mathcal{D}^-) = -1$. We call \emph{positive} the points labeled 1 and \emph{negative} the remaining ones. When we need a uniform naming convention for positive and negative points we use that $\boldsymbol{a}_i := \boldsymbol{a}^+_i$ for every $i \in I^+$ and $\boldsymbol{a}_i := \boldsymbol{a}^-_{i-m}$ for every $i \in I \setminus I^+$.

We say that a hyperplane separates the positive points from a certain subset of negative points, indexed by $\bar{I}^- \subseteq I^-$, if it admits a valid equation $(b, \boldsymbol{w})$ such that $b + \sum_{j=1}^d w_j a_{ij} < 0$ for every $i \in \bar{I}^-$. The \emph{indicator vector} of $\bar{I}^-$ is the vector $\boldsymbol{x} \in \{0, 1\}^n$ such that for each $i \in I^-$, $x_i = 1$ if and only if $i \in \bar{I}^-$.

A solution $s$ of \CPS separates the positive points from all the negative points if and only if $err(s) = 0$. The minimum hyperplane budget $\mathcal{K}$ such that \CPS admits a solution with separation error 0 is called \emph{minimum separation budget}.

The concatenation of two real-valued vectors $\boldsymbol x$ and $\boldsymbol y$ is indicated as $\boldsymbol x \mathbin\Vert \boldsymbol y$. Their element-wise product is written $\boldsymbol x \circ \boldsymbol y$. \\

\noindent
{\bf Related Work. }
Despite its simplicity, we are not aware of algorithmic approaches specifically targeting our \CPS. At the same time, it links to several relevant problems from the literature, not only from computational geometry, but also from machine learning and optimization. \\

Concerning computational geometry, our \CPS is related to the problem of computing the convex hull of a set of points in terms of the equations of the delimiting hyperplanes. It is well know in the literature: dedicated software packages were also conceived, engineered and maintained by the community \citep*{assarf2017computing}. Computing the convex hull in arbitrary dimension can be done exactly applying the algorithm QuickHull \citep*{Barber1996TheQA, Bykat1978ConvexHO, Eddy1977ANC} 
or approximately using, for example, more recent algorithms \citep*{Balestriero2022DeepHullFC}. Solving \CPS gives an over-approximation of the convex hull of the positive points but, on the other hand, \CPS focuses on finding much simpler separation boundaries by setting a hyperplane budget. It is worth noting, for instance, that while convex hull algorithms always produce polytopes, our \CPS models also allow unbounded convex polyhedra as valid solutions. \\

Another related problem is learning intersections of halfspaces using points labeled as positive or negative. As in the \CPS, the focus in on identifying the structure of a specific unknown polyhedron containing the positive points but none of the negative ones. However, unlike the \CPS, an implicit representation of the polyhedron is normally given, for example, by means of an oracle that augments the learning process. In this setting, polynomial-time algorithms w.r.t the dimension and the number of faces of the polyhedron were given by \cite{blum1995learning}, \cite{baum1990polynomial} and \cite{kwek1996pac}. Learning from random examples without access to an oracle, as in our case, is proved to be a challenging task whose related complexity results can be independent of the distribution from which the data is drawn \citep*{klivans2004learning} or distribution-dependent, as in previous works by \cite{baum1991neural} and \cite{vempala1997random}. \\

A close connection also exists with convex polyhedral classification (which is sometimes named convex polyhedral separation in the literature), a supervised learning task that consists in finding a binary classifier whose positive region is a convex polyhedron. It is a particular case of binary classification with piece-wise linear separation boundary. A foundational work is \citep*{Astorino2002PolyhedralST}, in which a weighted sum of the classification error on the positive class and on the negative one is minimized solving, at each iteration, a linear program providing a descent direction. The problem has received a fair amount of attention over the years, with solution approaches that cover mixed-integer programming \citep*{Orsenigo2007AccuratelyLF}, statistical learning \citep*{Manwani2010LearningPC} local SVM-like classification \citep*{Leng2020ASC}, and many others \citep*{Astorino2023MaximummarginPS, Astorino2021PolyhedralSV, Astorino2015SupportVM, Gottlieb2018LearningCP, Kantchelian2014LargeMarginCP, Raviv2018HingeMinimaxLF, Strekalovsky2015OnTP}. 
However, the focus of convex polyhedral classification is on finding a classifier with good generalization properties. Unlike our \CPS, the positive and negative sets are known only partially at training time. Therefore, they typically aim at separation boundaries with margins as large as possible, while \CPS tries to find convex hull approximations which are as tight as possible. Still some works in this area can be adapted to our case, as discussed in the remainder. \\

We finally mention the link to linear constraint learning. Some algorithms have been proposed that provide explicit representations of linear constraints based on examples of feasible and infeasible solutions, e.g.~\citep*{Pawlak2017AutomaticSO, Schede2019LearningLP}. These works represent just a sketch of the broader area of constraint learning. The interested reader can refer to \citep*{fajemisin2024optimization} for a general survey and the description of a conceptual framework for development of constraint learning techniques. 
Indeed, our \CPS algorithms are designed with a perspective application in a constraint learning framework. On one side as observed in \citep*{Messana2024Actively}, pure optimization algorithms as ours are hardly effective when used \emph{alone} for constraint separation, because there is no guarantee that available feasible and infeasible solutions are uniformly distributed. On the other side, when optimization is coupled with proper \emph{solutions sampling}, approaches like \CPS become promising.\\


\section{Models}
\label{sec:CompactModels}

We introduce two compact mixed-integer programming models for \CPS.
The first one is directly inspired by Support Vector Machines (SVMs), borrowing from earlier works on polyhedral separation. The second one incorporates key ideas to improve the approximation performance.

A classical mathematical programming model for the training of SVMs is reported in Appendix~\ref{sect:svm-model}. It uses one continuous variable $w_j$ for each dimension, encoding coefficients of a hyperplane, and one continuous (non-negative) variable $\xi_i$ to track the absolute misclassification error on point $i \in I$.

In fact, the SVM model can be adapted to solve \CPS as follows. We name this formulation as Model A, which is reported in full in Appendix~\ref{sect:svm-model}. Model A prescribes to find $\mathcal{K}$ hyperplane equations, encoded by variables $w_j^k$ with $j \in J$ and $k \in K$. A set of non-negative variables $\xi_i^k$ encode the distance from point $i \in I^-$ to hyperplane $k$, if $i$ is on the wrong side of the hyperplane $k$; they take value $0$ otherwise. Since \CPS requires every hyperplane in a solution to be valid, we neglect the violation variables corresponding to positive points, imposing them to be always on the correct side of each hyperplane.

We introduce binary variables $x_i^k$ for $i \in I^-$ and $k \in K$, having the following meaning: when $x_i^k = 1$, it means that point $i$ is \emph{separated} by the $k^{th}$ hyperplane. That is, we impose that $x_i^k$ can be set to $1$ only if the point $i \in I^-$ is on the negative side of hyperplane $k$. Each negative point must be separated by one hyperplane. We still allow negative points to be on the positive side of all hyperplanes (thereby being misplaced), at the price of a violation penalty. The objective function sums up the weighted violations over the negative points. For each negative point, the violation cost is considered only with respect to a single hyperplane (that yielding minimum violation). The weights are expressed as values of a function $\hat c: \mathbb{R}^d \to [0, +\infty)$.

The potential and limits of using binary variables in SVM-like models are discussed in the literature \citep*{carrizosa10}. Additionally, as we show in Section \ref{subsec:results}, the solutions to \CPS obtained solving Model A using Gurobi are weak in terms of approximation of the convex hull. It is not surprising, since the objective function of Model A targets a margin maximization, as in SVMs. Furthermore, it is quadratic, involving the product of variables $x_i^k$ and $\xi_i^k$. As final formulation step, therefore, we introduce the following mixed-integer programming model, denoted Model B.

\begin{align}
\text{min } & \sum_{i \in I^-} c(\boldsymbol{a}^-_i) e_i \label{obj:model_2} \\
\text{s.t. } & \sum_{k \in K} x_i^k = 1 - e_i & \forall i \in I^- \label{st:model_b_x_e} \\
            & b^k + \sum_{j \in J} w_j^k a_{ij} \geq 1 & \forall i \in I^+, \forall k \in K \label{st:model_b_positive} \\
            & b^k + \sum_{j \in J} w_j^k a_{ij} \leq - 1 + \xi_i^k & \forall i \in I^-, \forall k \in K \label{st:model_b_negative} \\
            & \xi_i^k \leq M \cdot (1-x_i^k) & \forall i \in I^-, \forall k \in K \label{st:model_b_big_M} \\
            & e_i \in \{0, 1\} & \forall i \in I^- \label{st:model_b_domain_e} \\
            & \xi_i^k\geq 0 & \forall i \in I^-, \forall k \in K \label{st:model_b_domain_xi} \\
            & x_i^k \in \{0, 1\} & \forall i \in I^-, \forall k \in K \label{st:model_b_domain_x} \\
            & b^k \in \mathbb{R} & \forall k \in K \label{st:model_b_domain_b} \\
            & w_j^k \in \mathbb{R} & \forall j \in J, \forall k \in K \label{st:model_b_domain_w}
\end{align}
Each binary variable $e_i$, with $i \in I^-$, represents the condition of negative point $i$ to be incorrectly included in the convex hull. When $e_i = 0$, constraints ~\eqref{st:model_b_x_e} impose that exactly one $x_i^k = 1$; in turn, when $x_i^k = 1$, constraints ~\eqref{st:model_b_big_M} impose that $\xi_i^k = 0$, and constraints ~\eqref{st:model_b_negative} impose $i$ to be on the correct side of hyperplane $k$. Constraints ~\eqref{st:model_b_positive} impose all positive points $i \in I^+$ to be on the correct side. The objective function \eqref{obj:model_2} counts the number of misplaced negative points (possibly weighting each of them differently). As technical remarks, the weights in the objective function do not need to be non-negative (while they do in Model A). Therefore, we can replace the weight function $\hat c$ with a real valued function $c: \mathbb{R}^d \to \mathbb{R}$. Furthermore, to enhance the effect of counting the misplaced points, saving one hyperparameter, no regularization term on $w$ was placed in the objective function. Accordingly, the right hand sides of constraints ~\eqref{st:model_b_positive} and \eqref{st:model_b_negative} contain constants ($1$ and $-1$) to avoid trivial solutions with $w = 0$. Finally, the big-M terms in constraints ~\eqref{st:model_b_big_M} have been introduced to avoid quadratic terms in the objective function.\\

In order to improve the bounds that can be obtained from Model B we employed Dantzig-Wolfe decomposition.  Let 
\[\mathcal{F}^k = \left\{\left(x_i^k\right)_{i \in I^-} \mathbin\Vert \left(w^k_j\right)_{j \in J} : \eqref{st:model_b_positive}, \eqref{st:model_b_negative}, \eqref{st:model_b_big_M}, \eqref{st:model_b_domain_xi}, \eqref{st:model_b_domain_x}, \eqref{st:model_b_domain_b}, \eqref{st:model_b_domain_w}\right\}\]
and let $\mathcal{G}^k = \left\{\left(x_i^k\right)_{i \in I^-}\ |\ \exists \left(w^k_j\right)_{j \in J}\ s.t.\ \left(x_i^k\right)_{i \in I^-} \mathbin\Vert \left(w^k_j\right)_{j \in J} \in \mathcal{F}^k\right\}$ for each $k \in K$.  Since $\mathcal{G}^k$ is a set of integer points, as follows:
\begin{align}
\text{min } & \sum_{i \in I^-} c(\boldsymbol{a}^-_i) e_i & \notag \\
\text{s.t. }& \sum_{k \in K} x_i^k = 1 - e_i & \forall i \in I^- \notag\\
            & \left(x_i^k\right)_{i \in I^-} \in conv\left\{\mathcal{G}^k\right\} & \forall k \in K \notag \\
            & e_i \in \{0,1\} & \forall i \in I^- \notag
\end{align}
We relax the integrality constraints on the error variables obtaining $0 \leq e_i \leq 1\ \forall i \in I^-$.

Let $\Omega^k$ be the set of extreme (integer) points of conv$\{\mathcal{G}^k\}$, and we denote (the coordinates of) each of them as $\bar x_i^{kp}$ for every $p \in \Omega^k$. Valid values for $x_i^k$ can be obtained as a convex combination of these extreme points. Let us introduce multipliers $z^{kp}$ representing such a convex combination, i.e. 
\begin{equation}
\label{st:convexdef}
x_i^k = \sum_{p \in \Omega^k} \bar x_i^{kp} z^{kp}\ \ \ \forall i \in I, \forall k \in K
\end{equation}
with $\sum_{p \in \Omega^k} z^{kp} = 1$ for each $k \in K$.
Constraints \eqref{st:convexdef} can be used to project out the $x_i^k$ variables:
\begin{align}
\text{min } & \sum_{i \in I^-} c(\boldsymbol{a}^-_i) e_i & \notag \\
\text{s.t. }& \sum_{k \in K} \sum_{p \in \Omega^k} \bar x_i^{kp} z^{kp} = 1 - e_i & \forall i \in I^-  \notag\\
            & \sum_{p \in \Omega^k} z^{kp} = 1 & \forall k \in K \notag \\
            & 0 \leq z^{kp} \leq 1 & \forall k \in K, \forall p \in \Omega^k \notag \\
            & 0 \leq e_i \leq 1 & \forall i \in I^- \notag
\end{align}
By observing the sets $\mathcal{G}^k$, one can notice all of them to contain the same set of points, call it $\mathcal{G}$ (and for consistency, denote as $\Omega$ the set of extreme points of the corresponding convex hull). Our model can be further simplified as follows:
%
%
\begin{align}
\text{min } & \sum_{i \in I^-} c(\boldsymbol{a}^-_i) e_i & \label{obj:rmm} \\
\text{s.t. }& \sum_{p \in \Omega} \bar x_{i}^p \cdot z^p \geq 1 - e_i & (\lambda_i)\ \  & \forall i \in I^-  \label{st:rmm_lambda} \\
            & -\sum_{p \in \Omega} z^p \geq -k & (\mu)\ \  & \label{st:rmm_mu} \\
            & 0 \leq z^p \leq 1 & & \forall p \in \Omega \label{st:rmm_domain_z} \\
            & 0 \leq e_i \leq 1 & & \forall i \in I^- \label{st:rmm_domain_e}
\end{align}
where $\lambda_i$ for all $i \in I^-$ and $\mu$ represent the dual variables associated to constraints \eqref{st:rmm_lambda} and \eqref{st:rmm_mu} respectively. Constraints \eqref{st:rmm_mu} are relaxed without loss of quality in the bound, to improve numerical stability.\\

\section{Algorithms}
\label{sect:algorithms}

The formulation \eqref{obj:rmm} - \eqref{st:rmm_domain_e} contains a combinatorial number of variables $z^p$. In the following, we build an algorithm that encompasses the standard column generation approach \citep*{amaldi12}, using such an extended formulation as a starting point. In addition, we propose an algorithm to solve the pricing that is specific to the problem, as it exploits its combinatorial structure. 
The reduced cost of a variable $z^p$ can be written as $\bar{z}^ p = 0 - \sum_{i \in I^-} \bar x_i^p \cdot \lambda_i + \mu$. The pricing model (PM) is therefore the following:
\begin{align}
\text{min } \bar{z}^ p =  & - \sum_{i \in I^-} x_i \cdot \lambda_i + \mu \label{obj:PM} \\
\text{s.t. } & b + \sum_{j \in J} w_j a_{ij} \geq 1 & \forall i \in I^+ \\
            & b + \sum_{j \in J} w_j a_{ij} \leq - 1 + \xi_i & \forall i \in I^- \\
            & \xi_i \leq M \cdot (1-x_i) & \forall i \in I^- \\
            & \xi_i \in [0, +\infty) & \forall i \in I^- \label{st:domain_xi} \\
            & x_{i} \in \{0,1\} & \forall i \in I^-\\
            & b \in \mathbb{R} \\
            & w_j \in \mathbb{R} & \forall j \in J \label{st:wPM}
\end{align}
We explicitly indicate the dependence on the specific dual variables by writing PM($\boldsymbol{\lambda}$, $\mu)$. \\

A restriction of model \eqref{obj:rmm} - \eqref{st:rmm_domain_e} obtained by replacing $\Omega$ with a subset $\Omega^{\prime} \subseteq \Omega$, is called Restricted Master Model (RMM). Given a RMM, we call Hyperplane Choice Model (HCM) the same model with 0-1 variables instead of continuous ones.

Our version of column generation with hyperplane choice is detailed in Alg.~\ref{alg:CGandHC} of Appendix \ref{sec:ColumnGenerationAndHyperplaneChoice}. At each iteration, we solve the restricted master model and the pricing problem. We address the pricing problem applying one of two different approaches, solving PM with general purpose solvers, or applying the algorithm that we show later in Section \ref{subsec:ahp}. Subsequently, RMM is updated adding the columns obtained from the solution of the pricing problem and HCP is solved to find the best possible solution to \CPS using hyperplanes obtained from solving pricing. This last step allows us to monitor the error value at each iteration. The algorithm returns with the last available \CPS solution when all the columns generated in a certain iteration have non-negative reduced cost or when the error value of the \CPS solution is 0.

Notice that, when solving PM using general purpose solvers, we get by default a single column. However, Alg.~\ref{alg:CGandHC} takes into account the possibility to employ the algorithm AHP of Section \ref{subsec:ahp} instead, which can generate more than one column at each iteration.

\subsection{Ad Hoc Heuristic Algorithm for the Pricing Problem}
\label{subsec:ahp}

Solving the pricing problem implies placing a single valid hyperplane. For the solution to be optimal, we want to maximize the sum of the dual variables associated to the negative points that are separated from the positive points by means of that hyperplane. A subset of negative points is separable from the set of positive points if and only if their respective convex hulls are disjoint. Therefore, in order to find a good heuristic solution to the pricing problem, it is sufficient to find a subset of negative points such that:

\begin{enumerate}
    \item the sum of the corresponding dual variables is ``high'';
    \item its convex hull does not intersect the convex hull of the positive points.
\end{enumerate}

The main intuition behind our ad hoc heuristic algorithm for the pricing problem is that such a set of negative points can be found in a sequential fashion, considering the negative points one at a time and keeping only those that maintain the second condition true. We try to ensure the first condition considering the negative points accordingly to their dual variable values, in decreasing order. This approach is appealing in that it lends itself to parallelization. In our actual implementation, we run twin parallel processes of Alg.~\ref{alg:hpls} by changing the starting point only; this increases the chances of obtaining a good solution. \\


Given a subset $\bar{I}^- \subseteq I^-$ of negative points, it is possible to check if the conxex hull of the points indexed by $\bar{I}^-$ intersects the convex hull of the positive points after solving the following linear programming model, called linear separability verification model and denoted LSVM($\bar{I}^-$). \\

\begin{align}
\text{min } & \sum_{j \in J} \left(\xi_j + \eta_j\right) \label{obj:LSV} \\
\text{s.t. } & \sum_{i \in I^+} \alpha_i \cdot a^+_{ij} + \xi_j = \sum_{i \in \bar{I}^-} \beta_i \cdot a^-_{ij} + \eta_j & \forall j \in J \\
             & \sum_{i \in I^+} \alpha_i = 1 \\
             & \sum_{i \in \bar{I}^-} \beta_i = 1 \\
             & \alpha_i \in [0,\ 1] & \forall i \in I^+ \\
             & \beta_i \in [0,\ 1] & \forall i \in \bar{I}^- \\
             & \xi_j \geq 0 & \forall j \in J \\
             & \eta_j \geq 0 & \forall j \in J \label{st:etaLSV}
\end{align}
The model admits a solution with all zero slacks, i.e. $\xi_j = 0$ and $\eta_j = 0$ for every $j \in J$, if and only if it exists a convex combination of the positive points equal to a convex combination of the given negative points or, in other words, if their respective convex hulls intersect. \\

The pseudocode of our ad hoc algorithm is given in Alg.~\ref{alg:hpls}
 and Alg.~\ref{alg:ahp}. Alg.~\ref{alg:hpls} shows the process called \emph{heuristic partial linear separation} (HPLS), whose role is to find the desired set of negative points and to derive a single pricing solution, given by an indicator vector $\boldsymbol x$ and hyperplane parameters $b$ and $\boldsymbol{w}$. It requires as input the index $i^{start} \in I^-$ of a starting negative point, the maximum number $thr$ of negative points to consider and a permutation $\sigma: I^- \rightarrow I^-$ such that $\lambda_{\sigma(i)} \geq \lambda_{\sigma(i')}$ for every $i, i' \in I^-$ with $i < i'$.
\begin{algorithm}
\small
\caption{Heuristic Partial Linear Separation}\label{alg:hpls}
\begin{algorithmic}
    \State \textbf{Input:} permutation $\sigma$, starting index $i^{start}$, point threshold $thr$
    \State
    \State $\bar{I}^- \gets \{i^{start}\}$
    \State $counter \gets 0$
    \While{$counter < n$ \textbf{and} $counter < thr$}
        \If{$\sigma(counter) = i^{start}$}
            \State $counter \gets counter + 1$
            \State \textbf{continue}
        \EndIf
        \State Solve LSVM($\bar{I}^-$) and get slack variable values $\xi_j$ and $\eta_j$ for every $j \in J$
        \If{$\exists j \in J$ such that $\left(\xi_j > 0\ \textbf{or}\ \eta_j > 0\right)$}
            \State $\bar{I}^- \gets \bar{I}^- \cup \{\sigma(counter)\}$
        \EndIf
    \EndWhile
    \State Solve a linear separation problem between all the positive points and the negative points indexed by $\bar{I}^-$ and get hyperplane coefficients $b$ and $\boldsymbol{w}$
    \State Shift the hyperplane as close as possible to the set of positive points and get new coefficient b
    \State $\boldsymbol{x} \gets \boldsymbol{0} \in \mathbb{R}^{n}$
    \For{$i \in I^-$}
        \If{$b + \langle \boldsymbol{w},\ \boldsymbol{a}^-_i \rangle \ \leq -1$}
            \State $x_i \gets 1$
        \EndIf
    \EndFor \\
    \Return $\boldsymbol{x}$, $b$, $\boldsymbol{w}$
\end{algorithmic}
\end{algorithm}
Until all the negative points have been considered or the point threshold $thr$ has been reached, a new point is added to $\bar I^-$ and the model LSVM($\bar{I}^-$) is solved to determine whether to keep or to neglect the point. After the end of the loop, we solve a linear separation problem to find an equation $(b, \boldsymbol{w})$ of a hyperplane that separates the negative points indexed by the final $\bar{I}^-$ from the positive points. This hyperplane is then moved as close as possible to the set of positive points, in order to possibly increase the number of negative points that it classifies correctly. The indicator vector $x$ is obtained assigning 1 to every negative point in $\bar{I}^-$ and 0 to the other ones.

\begin{algorithm}
\small
\caption{Ad Hoc Heuristic Algorithm for the Pricing Problem}\label{alg:ahp}
\begin{algorithmic}
    \State \textbf{Input:} dual variables $\boldsymbol{\lambda} = (\lambda_i)_{i \in I^-}$ and $\mu$, point threshold $thr$, vector $\boldsymbol{x}^{prev} \in \{0, 1\}^{n}$, maximum number of processes $N^{max}$
    \State
    \State Find a permutation $\sigma: I^- \rightarrow I^-$ such that $\lambda_{\sigma(i)} \geq \lambda_{\sigma(i')}$ for every $i, i' \in I^-$ with $i < i'$
    \State Set $\tilde x_i = 0$ for every $i \in I^-$
    \State $\boldsymbol{p} \gets \boldsymbol{x}^{prev} \circ \boldsymbol{\lambda}$
    \State $p_i \gets p_i / \sum_{j \in I^-} p_j$ for every $i \in I^-$
    \State Draw $i^{start}$ from $I^-$ according to distribution $\boldsymbol{p}$
    \If{$\exists i \in I^-$ such that $\bar{x}_i^{prev} = 0$}
        \State $\bar{y}_i^{prev} \gets 1 - \bar{x}_i^{prev}$ for every $i \in I^-$
    \Else
        \State $\bar{y}_i^{prev} \gets 1$ for every $i \in I^-$
    \EndIf
    \State $\boldsymbol{p} \gets \bar{\boldsymbol{y}}^{prev} \circ \boldsymbol{\lambda}$
    \State $p_i \gets p_i / \sum_{j \in I^-} p_j$ for every $i \in I^-$
    \State $N \gets 0$
    \While{$N < N^{max}$ \textbf{and} $\exists i \in I^-$ such that $\tilde x_i = 0$}
        \State Run parallel HPLS($\sigma$, $i^{start}$, $thr$) process
        \State $N \gets N + 1$
        \If{$N = N^{max}$}
            \State \textbf{break}
        \EndIf
        \State Draw $i^{start}$ from $I^-$ according to distribution $\boldsymbol{p}$
    \EndWhile
    \State Wait for the running HPLS processes to terminate and collect $\boldsymbol{x}^k$, $\boldsymbol{w}^k$ and $b^k$ for every process $k \in \bar{N}$ \\
    \Return $N$ and $\boldsymbol{x}^k$, $b^k$, $\boldsymbol{w}^k$ for every $k \in \bar{N}$
\end{algorithmic}
\end{algorithm}

Alg.~\ref{alg:ahp} shows the general structure of the ad hoc heuristic for pricing (AHP), whose task is to obtain multiple pricing solutions by running parallel HPLS processes. The following input parameters are required: $\boldsymbol{\lambda}$ and $\mu$ are the dual variables obtained from the solution of the most recent restricted master model, $thr$ is the point threshold to be used in HPLS, $\boldsymbol{x}^{prev}$ is the indicator vector corresponding to the pricing solution with lowest reduced cost found in the previous iteration of column generation and finally, $N^{max}$ is the maximum allowed number of parallel instances of HPLS. Notice that, according to Alg.~\ref{alg:CGandHC}, at the first iteration of column generation, $\boldsymbol{x}^{prev} = \boldsymbol{1}$. \\

As already mentioned, every HPLS process requires a starting index $i^{start}$. Such starting indices are chosen at random with probability proportional to the corresponding values of the dual variables $\boldsymbol{\lambda}$, in order to favour points with higher dual variable values. For the first parallel process, only indices with a corresponding 1 in the indicator vector $\boldsymbol{x}^{prev}$ can be chosen. For the subsequent processes, it holds the opposite: an index can be chosen only if the corresponding value in $\boldsymbol{x}^{prev}$ is 0. In other words, for the first process we select the index of a point that, in the previous iteration of column generation, has been separated from the positive points by the best hyperplane (one with lowest corresponding reduced cost among the ones obtained); for the other processes, we choose only indices of negative points that were erroneously classified as positive by the same hyperplane. This differentiated behaviour allows us to try an improvement mimicking the best pricing solution of the previous iteration and, at the same time, to put effort in exploring different solutions.

\subsection{A Greedy Matheuritic for \CPS}
\label{sec:GreedyAlgorithm}

We devised also a greedy matheuristic approach to solve \CPS, for both obtaining an initial solution quickly, and to populate the initial RMP. It consists in placing a single hyperplane at a time and ignoring every negative point that is already separated correctly by at least one of the hyperplanes created so far. To find a hyperplane that, given a certain subset $\bar{I}^- \subseteq I^-$ of negative points, separates the positive points from the greatest number of negative points in $\bar{I}^-$, it is possible to solve the model below, called \emph{Partial Linear Separation Model}.
\begin{align}
\text{max} & \sum_{i \in \bar{I}^-} c(\boldsymbol{a}^-_i)x_i \\
\text{s.t. } & b + \sum_{j \in J} w_j a_{ij} \geq 1 & \forall i \in I^+ \\
            & b + \sum_{j \in J} w_j a_{ij} \leq - 1 + \xi_i & \forall i \in \tilde I^- \\
            & \xi_i \leq M \cdot (1-x_i) & \forall i \in \tilde I^- \\
            & \xi_i \in [0, +\infty) & \forall i \in \tilde I^- \\
            & x_i \in \{0, 1\} & \forall i \in \tilde I^- \\
            & b \in \mathbb{R} \\
            & w_j \in \mathbb{R} & \forall j \in J
\end{align}
We write PLSM($\bar{I}^-$) to denote the model and its dependence on the set $\bar{I}^-$. This model is similar to \eqref{obj:PM} - \eqref{st:wPM}, the pricing model of our column generation. The differences are that not necessarily all the negative points are involved and that the objective is a particular case of \eqref{obj:PM}, in which $\lambda_i = 1$ for all $i$ and $\mu = 0$. \\

In Alg.~\ref{alg:greedy} we show a pseudocode for the greedy algorithm.
\begin{algorithm}
\small
\caption{Greedy Algorithm}\label{alg:greedy}
\begin{algorithmic}
\State Set $x_i = 0$ for every $i \in I^-$
\State $\bar I^- \gets I^-$
\State $t \gets 0$
\While{$t < \mathcal{K}$ \textbf{and} $\exists i \in I^-$ such that $x_i = 0$}
    \State Solve PLSP($\bar I^-$) and get $\{x_i^t\}_{i \in \bar I^-}$, $b^t$ and $\boldsymbol{w}^t = (w_j^t)_{j \in J}$
    \State Set $x_i = 1$ for every $i \in \bar I^-$ such that $x_i^t = 1$
    \State $\bar I^- \gets \{i \in \bar I^-\ |\ x_i = 0\}$
    \State $t \gets t + 1$
\EndWhile \\
\Return $\{(\boldsymbol{w}^{\tau},\ b^{\tau})\ |\ \tau \in \{1, \ldots, t\}\}$
\end{algorithmic}
\end{algorithm}
Binary variables $x_i$ indicate whether the negative point $i$ is classified correctly by at least one chosen hyperplane or not. The set $\bar{I}^-$ is initialized to contain all the negative points. As long as less that $\mathcal{K}$ hyperplanes have been placed and there is still some misclassified negative points, a hyperplane is chosen solving PLSM($\bar{I}^-$) and all the points classified correctly by that hyperplane are removed from $\bar{I}^-$.

\section{Experimental Setting and Results}
\label{sec:ExperimentalSettingAndResults}

 We create two datasets (D1 and D2), each composed by 3 sets of 100 instances each, with number of dimensions $d$ equal to 2, 4 and 8 respectively, assuming the unitary hypercube as the feasible region and generating an exponential number of points with respect to the dimension (ranging from about 200 points for $d = 2$ to 10000 points for $d = 8$). D1 and D2 are generated similarly, but instances in D1 have a known minimum separation budget, whereas instances in D2 do not. Details on the generations of our instances are reported in Appendix \ref{sec:DataGeneration}; both datasets are available online \citep*{RD_UNIMI/IE2MX0_2024}.
 The experimental analysis is carried out in two steps. In the first part, we use dataset D1 and we consider our compact models \eqref{obj:model_a} - \eqref{st:model_a_domain_w} (named MODEL A) and \eqref{obj:model_2} - \eqref{st:model_b_domain_w} (MODEL B), column generation with original pricing formulation (COLGEN 1) and with ad hoc pricing algorithm (COLGEN 2), and finally the greedy algorithm of Section \ref{sec:GreedyAlgorithm} (GREEDY). We also adapted a compact model proposed in \citep*{Orsenigo2007AccuratelyLF} for convex polyhedral classification; the original model was designed for a problem closely matching ours, and was therefore easy to adapt. Such an adaptation, that we name OV 2007, is detailed in Appendix \ref{sec:AppendixOV2007}.

In the second part of the analysis, we investigate the capability of a column generation approach to approximate the convex hull of the positive points, using dataset D2. In particular we consider COLGEN 2 and we compare it against a procedure customized from the state-of-the art software polymake, \cite{assarf2017computing}. For this second part of the experiments, we use very similar datasets as in the first part, but with unknown minimum separation budget.

All the programming models, including the ones involved in column generation and in the greedy algorithm, have been tackled using Gurobi. The code was written and executed in Python 3.9 on a machine with Ubuntu 20.04.6 LTS, 32 GB DDR4 RAM and processor Intel(R) Xeon(R) W-1250P CPU @ 4.10GHz.

\subsection{Implementation Choices}

For Model B, we fixed $M = 10000$, which proved to be a sufficiently large value not to interfere negatively with the solution process. For OV 2007, we tried different values for $\epsilon$, including machine precision as suggested by the authors in the original paper. However, despite several attempts, we have not been able to find a value that always ensures a correct behaviour in the solution process and in the end we set $\epsilon = 0.001$. For further details, see Section \ref{subsec:results}.

Both the column generation algorithms considered are heuristic. The first one is Alg.~\ref{alg:CGandHC} with PM($\boldsymbol{\lambda}$, $\mu$) solved using Gurobi and interrupting the solver inside the branch-and-bound algorithm immediately after the root node. The reason for this choice is twofold: preliminary tests showed us that the solution process of PM terminates within a reasonable time only on very small datasets; moreover, the value of the objective function at the root node is often close to the values obtained going on in the branch-and-bound. The second column generation algorithm has the same structure of the first one, but AHP (Alg.~\ref{alg:ahp}) is executed instead of solving PM. We have chosen the maximum number of parallel processes $N^{max}$ equal to 8. The point threshold $thr$ has been set equal to $d$.

In the greedy algorithm, like in the first column generation, the solution process of PLSM is interrupted right after the root node of the branch-and-bound tree for the same reasons.

For all the methods, the weight function (that is, $\hat c$ for Model A and $c$ for the others) has been kept constant and equal to 1 on every $a \in \mathbb{R}^d$.



\subsection{Comparing Models and Algorithms}
\label{subsec:results}

We carried out preliminary experiments using 8 out of 100 instances for each of the three dimensions using large time limits, in order to have feedback about the performance of the models and the algorithms tested on the long run. Subsequently, we ran tests using all the 100 instances with lower time limits.

\paragraph{Preliminary experiments}

For each $d \in \{2, 4, 8\}$, we tested the six methods using three different hyperplane budgets: $d$, $2d - 1$ and $2d$. As already mentioned in Section \ref{sec:DataGeneration}, $\mathcal{K} = 2d$ represents the minimum separation budget. Each method has been run with a time limit equal to $6 \cdot 10^{log_2(d)}$ seconds (i.e.~$10^{log_2(d) - 1}$ minutes).

\begin{figure}
    \centering
    \includegraphics[width=0.7\textwidth]{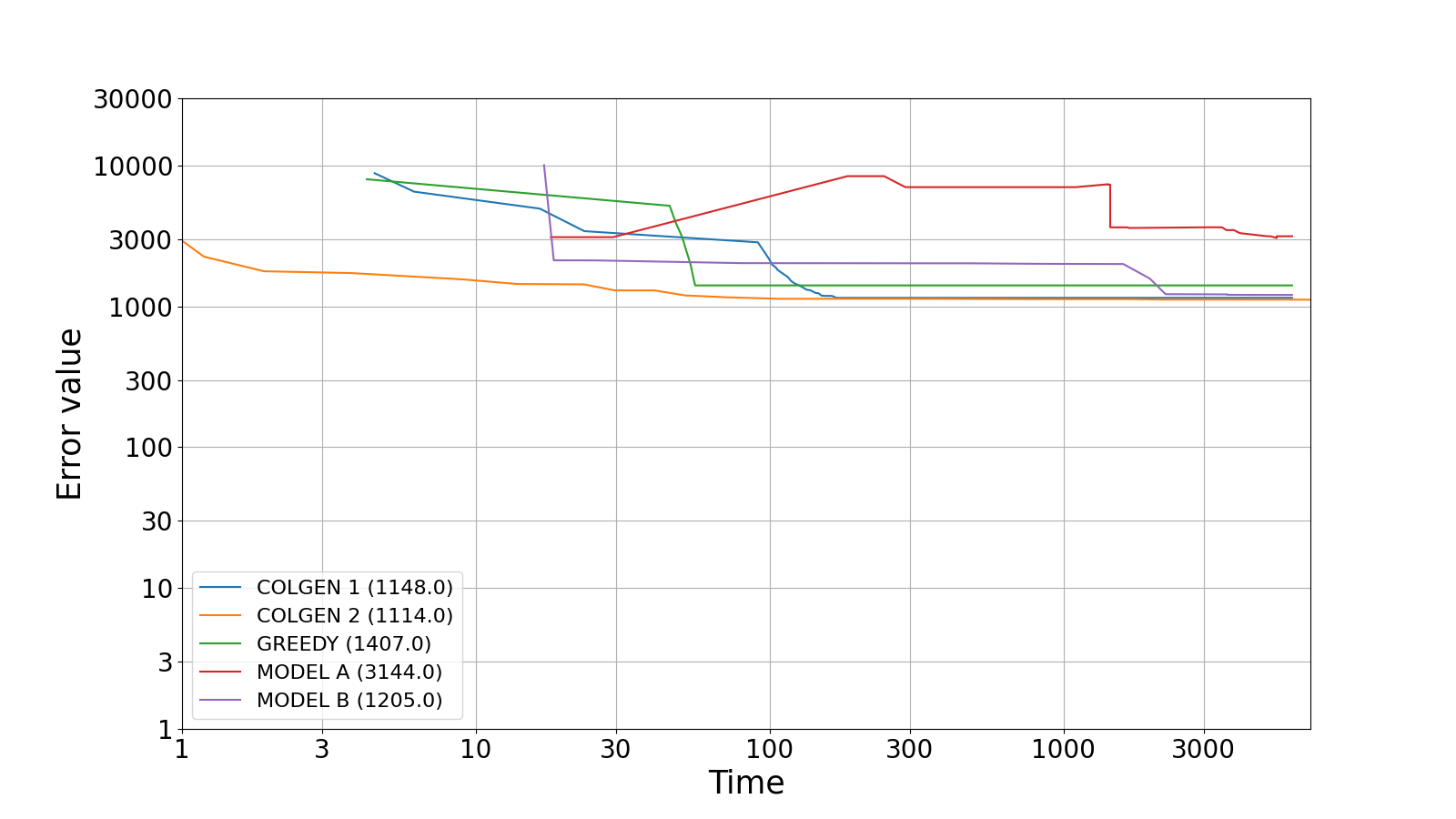}
    \caption{Error value with 8-dimensional hypercubic data: $\mathcal{K} = 8$ (sample instance).}
    \label{fig:d8K8}
\end{figure}

\begin{figure}
    \centering
    \includegraphics[width=0.7\textwidth]{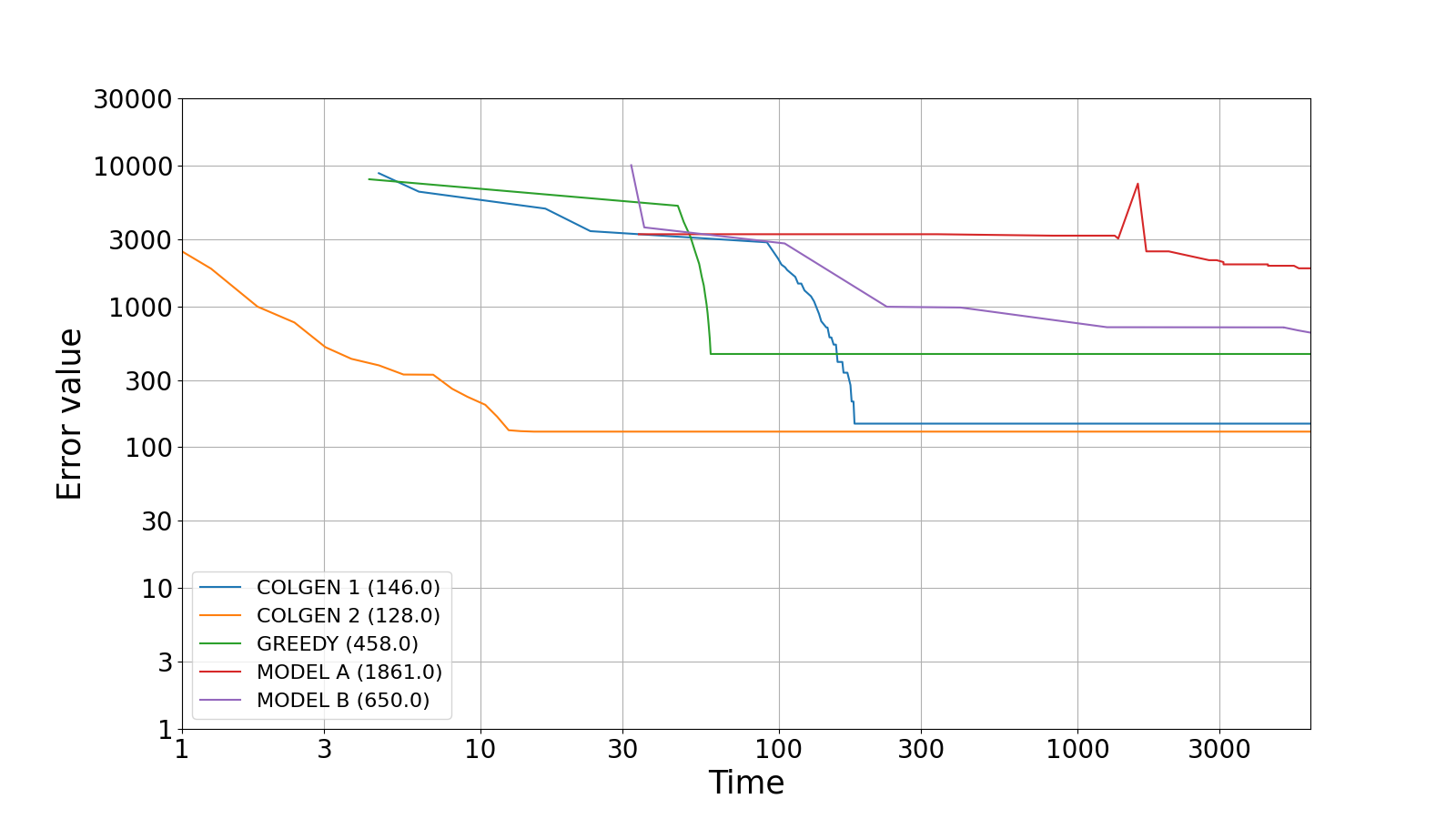}
    \caption{Error value with 8-dimensional hypercubic data: $\mathcal{K} = 15$ (sample instance).}
    \label{fig:d8K15}
\end{figure}

\begin{figure}
    \centering
    \includegraphics[width=0.7\textwidth]{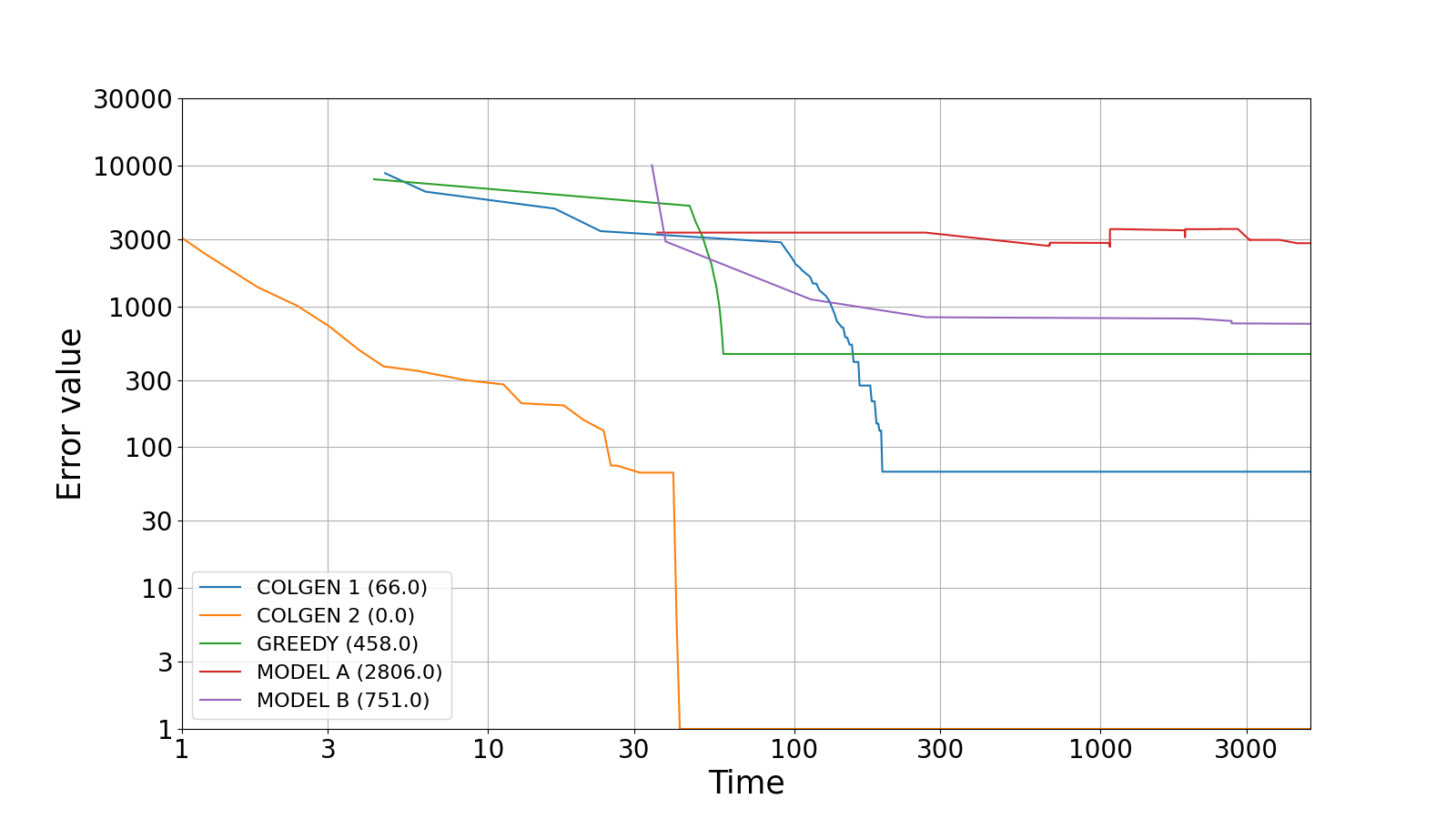}
    \caption{Error value with 8-dimensional hypercubic data: $\mathcal{K} = 16$ (sample instance).}
    \label{fig:d8K16}
\end{figure}

In order to give a representative example, Figures \ref{fig:d8K8} - \ref{fig:d8K16} show the error value over time for a single instance in dimension 8. Analogous figures are present in Appendix \ref{sec:AdditionalFiguresAndTables} for dimensions 2 and 4. Both error and time are reported in log-10 scale. In the legends, next to the method names, their respective final error values are indicated in brackets. In the following, we state the main observations that we can derive from the figures, but very similar comments hold for all the instances in the preliminary experiments, according to their dimensionality. \\

First of all, we can notice that the performance of Model A in terms of error value is the worst among all the methods in most of the cases the unique exception is the experiment whose results are reported Fig.~4). So we can affirm that generalizing in a straightforward way a popular and vastly applied model such as SVM does not yield good results for our problem. Also notice that the lines for Model A are non-monotone, because the objective that it tries to minimize is different from the error value. The column generation algorithms yield the best performances overall. This holds especially for COLGEN 2, which converges very quickly and obtains the best error values in the majority of the cases. Model B and GREEDY give intermediate results. \\

For the instance in dimension 2, all the methods converge really quickly. In dimension 4 and 8, we can notice more differentiated convergence speeds. We can also see that the extent of the difference range for error value and convergence speed seems to be directly proportional to the hyperplane budget. Finally, the reader may notice that OV 2007 is not present in Figures \ref{fig:d8K8} - \ref{fig:d8K16}. The reason is that the solutions obtained from solving that model are not valid: the solver outputs $b = 0$ and $\boldsymbol{w} = \boldsymbol{0}$ as a computational consequence of the fact that we were not able to find a value of the hyperparameter $\epsilon$ that allows the solver to find a proper solution when $d = 8$.

\paragraph{Massive experiments}

In the preliminary experiments we were able to verify that one tenth of the given time limit was always sufficient for the column generation algorithms to reach approximately the final error value. Therefore, we tested the six methods on all the 100 instances using the same hyperplane budgets as in the preliminary experiments, but time limit equal to $6 \cdot 10^{log_2(d) - 1}$ for each $d \in \{2,\ 4,\ 6\}$. Tables \ref{tab:massive_2} - \ref{tab:massive_8} show the results averaged over all the instances. We report running time, final error value and time needed to reach such error value (this latter abbreviatd to TE, time up to error), for each of the three dimensions respectively. Additional tables can be found in Appendix \ref{sec:AdditionalFiguresAndTables}. They contain details about the performance of the two column generation algorithms.

\begin{table}[bt]
\small
\centering
\begin{tabular}{c l r r r}
$\mathcal{K}$ & method & time (seconds) & \% error value & TE (seconds) \\
\hline
\multirow{6}{0.2cm}{2} & COLGEN 1 & 6.044 & 22.28 & 0.625 \\
& COLGEN 2 & 6.052 & 22.66 & 2.747 \\
& GREEDY & \textbf{0.164} & 27.23 & \textbf{0.164} \\
& MODEL A & 6.072 & 33.87 & 3.843 \\
& MODEL B & 6.081 & \textbf{21.24} & 0.856 \\
& OV 2007 & 6.190 & 22.47 & 4.253 \\
\hline
\multirow{6}{0.2cm}{3} & COLGEN 1 & 6.040 & 7.87 & 0.767 \\
& COLGEN 2 & 6.052 & 6.52 & 3.304 \\
& GREEDY & \textbf{0.199} & 9.86 & \textbf{0.198} \\
& MODEL A & 6.102 & 13.36 & 3.876 \\
& MODEL B & 6.109 & \textbf{6.07} & 1.774 \\
& OV 2007 & 6.243 & 7.40 & 3.358 \\
\hline
\multirow{6}{0.2cm}{4} & COLGEN 1 & 1.408 & 0.11 & 0.711 \\
& COLGEN 2 & 0.474 & \textbf{0.00} & 0.472 \\
& GREEDY & \textbf{0.208} & 0.92 & \textbf{0.207} \\
& MODEL A & 5.500 & 1.52 & 4.179 \\
& MODEL B & 0.220 & \textbf{0.00} & 0.215 \\
& OV 2007 & 1.460 & 0.08 & 1.388 \\
\hline
\end{tabular}
\caption{Average running time, error value and final error hitting time for $d = 2$.}
\label{tab:massive_2}
\end{table}

\begin{table}[tb]
\small
\centering
\begin{tabular}{c l r r r}
$\mathcal{K}$ & method & time (seconds) & \% error value & TE (seconds) \\
\hline
\multirow{6}{0.2cm}{4} & COLGEN 1 & 60.108 & 12.48 & 8.409 \\
& COLGEN 2 & 60.457 & \textbf{12.07} & 29.835 \\
& GREEDY & \textbf{2.111} & 15.89 & \textbf{2.108} \\
& MODEL A & 60.410 & 24.66 & 52.789 \\
& MODEL B & 60.408 & 13.49 & 25.264 \\
& OV 2007 & 60.911 & 22.37 & 51.594 \\
\hline
\multirow{6}{0.2cm}{7} & COLGEN 1 & 58.990 & 1.96 & 7.568 \\
& COLGEN 2 & 60.142 & \textbf{1.44} & 4.159 \\
& GREEDY & \textbf{2.232} & 2.81 & \textbf{2.227} \\
& MODEL A & 60.702 & 13.17 & 57.571 \\
& MODEL B & 60.686 & 1.76 & 14.369 \\
& OV 2007 & 61.501 & 11.14 & 49.668 \\
\hline
\multirow{6}{0.2cm}{8} & COLGEN 1 & 32.357 & 0.36 & 7.760 \\
& COLGEN 2 & \textbf{1.185} & \textbf{0.00} & \textbf{1.180} \\
& GREEDY & 2.267 & 1.22 & 2.260 \\
& MODEL A & 60.781 & 15.92 & 55.100 \\
& MODEL B & 18.662 & 0.10 & 14.873 \\
& OV 2007 & 61.696 & 8.43 & 52.282 \\
\hline
\end{tabular}
\caption{Average running time, error value and final error hitting time for $d = 4$.}
\label{tab:massive_4}
\end{table}

\begin{table}[tb]
\small
\centering
\begin{tabular}{c l r r r}
$\mathcal{K}$ & method & time (seconds) & \% error value & TE (seconds) \\
\hline
\multirow{6}{0.2cm}{8} & COLGEN 1 & 601.465 & 11.20 & 95.704 \\
& COLGEN 2 & 746.987 & \textbf{11.08} & 517.667 \\
& GREEDY & \textbf{30.880} & 12.80 & \textbf{30.795} \\
& MODEL A & 618.026 & 60.33 & 442.587 \\
& MODEL B & 617.558 & 17.68 & 581.188 \\
& OV 2007 & 637.087 & - & 35.709 \\
\hline
\multirow{6}{0.2cm}{15} & COLGEN 1 & 433.135 & 2.87 & 108.615 \\
& COLGEN 2 & 603.572 & \textbf{1.27} & \textbf{31.553} \\
& GREEDY & \textbf{33.963} & 5.22 & 33.487 \\
& MODEL A & 633.393 & 33.83 & 46.513 \\
& MODEL B & 632.671 & 10.25 & 455.574 \\
& OV 2007 & 666.978 & - & 65.929 \\
\hline
\multirow{6}{0.2cm}{16} & COLGEN 1 & 288.734 & 2.66 & 109.489 \\
& COLGEN 2 & \textbf{26.953} & $<$\textbf{ 10}$^{-4}$ & \textbf{26.717} \\
& GREEDY & 34.245 & 4.78 & 33.687 \\
& MODEL A & 635.822 & 30.20 & 337.217 \\
& MODEL B & 634.875 & 10.20 & 506.273 \\
& OV 2007 & 671.311 & - & 70.104 \\
\hline
\end{tabular}
\caption{Average running time, error value and final error hitting time for $d = 8$.}
\label{tab:massive_8}
\end{table}

The massive experiments confirm the relative performances showed by the methods in the preliminary experiments. In dimension 2, the small size of the instances allows MODEL B and OV 2007 to reach error values lower or comparable to those obtained by column generation. In dimension 4, MODEL B, COLGEN 1 and COLGEN 2 all give small error values but we can notice that the time up to error can be significantly higher for MODEL B. Finally, in dimension 8, COLGEN 2 proves to be the method with best error values, while COLGEN 1 is slightly less performing. In decreasing order of final error, the other methods are GREEDY, MODEL B, MODEL A and OV 2007. Notice that greedy is the only method that always terminates within the given time limit, thanks to the stopping criterion of its iteration step. For the other methods, reaching optimality within the time limit turns out to be difficult except when $d$ is not too large and $\mathcal K=2d$. However, both column generation algorithms terminate within the time limit when $\mathcal{K} = 2d$ for every $d \in \{2, 4, 8\}$.

Summarizing, the experiments highlight COLGEN 2 as the most promising method overall. However, it is at the same time the most complex. Parallelization showed to be an essential component. Its effective use require the fine tuning of the point threshold \emph{thr}, which might not be easy in all applications.

\subsection{Producing Polyhedral Approximations}

In a second round of experiments we assessed the effectiveness of our algorithms in producing good convex hull approximations. We considered COLGEN 2, being the most effective algorithm in previous experiments, with a time limit set as in the previous section. As a benchmark we considered a procedure adapted from Polymake. Polymake implements convex hull algorithms, possibly producing way more than $k$ hyperplanes. We considered a greedy approach (GREEDY POLYMAKE): firstly, it computes the equations of the hyperplanes describing the facets of the convex hull of the positive points; secondly, it iteratively selects a hyperplane that correctly classifies the highest number of negative points and it removes those points from the dataset. An additional hyperplane gets selected with the same criterion, until the budget $\mathcal{K}$ is reached. A larger time limit of one hour was given to GREEDY POLYMAKE.

We compared COLGEN 2 and GREEDY POLYMAKE on all the 100 instances generated for this second part of the experiments, applying the same hyperplane budget as in part one. The metrics for the comparison are the error, i.e. the number of misclassified negative point, and the volume of the convex region delimited by the resulting hyperplanes. In the cases in which the hyperplanes generated by an algorithm for a specific instance describe an unbounded polyhedron, we consider the region bounds for the generated points, i.e. $-1 \leq x_i \leq 2$ for every $i \in \{1, \ldots, d\}$. The computation of the convex hulls for GREEDY POLYMAKE and the computation of the volume of the resulting polyhedron for both the methods were achieved using ake \citep*{gawrilow2000polymake, assarf2017computing}. Tables \ref{tab:p2d2}, \ref{tab:p2d4} and \ref{tab:p2d8} show the results in dimension 2, 4 and 8 respectively, averaged over the 100 instances for each hyperplane budget and each algorithm. Table \ref{tab:p2d8} doesn't contain those results that require to use ake, as GREEDY POLYMAKE was always hitting the time limit (without producing any valid solution) also computing the volume was extremely computational expensive in dimension 8.

\begin{table}[bt]
\small
\centering
\begin{tabular}{c l r r}
$\mathcal{K}$ & method & \makecell{\% error value} & \makecell{volume} \\
\hline
\multirow{2}{0.2cm}{2} & COLGEN 2 & 22.01 & 3.00 \\
& GREEDY POLYMAKE & 21.52 & 2.96 \\
\hline
\multirow{2}{0.2cm}{3} & COLGEN 2 & 5.30 & 1.69 \\
& GREEDY POLYMAKE & 9.44 & 1.91 \\
\hline
\multirow{2}{0.2cm}{4} & COLGEN 2 & 0.00 & 1.00 \\
& GREEDY POLYMAKE & 0.06 & 0.98 \\
\hline
\end{tabular}
\caption{Average error value and volume of the approximate convex hull for $d = 2$.}
\label{tab:p2d2}
\end{table}

\begin{table}[bt]
\small
\centering
\begin{tabular}{c l r r}
$\mathcal{K}$ & method & \makecell{\% error value} & \makecell{volume} \\
\hline
\multirow{2}{0.2cm}{4} & COLGEN 2 & 5.33 & 7.64 \\
& GREEDY POLYMAKE & 6.93 & 9.00 \\
\hline
\multirow{2}{0.2cm}{7} & COLGEN 2 & 0.05 & 2.11 \\
& GREEDY POLYMAKE & 0.35 & 1.98 \\
\hline
\multirow{2}{0.2cm}{8} & COLGEN 2 & 0.00 & 1.71 \\
& GREEDY POLYMAKE & 0.04 & 1.37 \\
\hline
\end{tabular}
\caption{Average error value and volume of the approximate convex hull for $d = 4$.}
\label{tab:p2d4}
\end{table}

\begin{table}[h!]
\small
\centering
\begin{tabular}{c l r r}
$\mathcal{K}$ & method & \makecell{\% error value} & \makecell{volume} \\
\hline
\multirow{2}{0.2cm}{8} & COLGEN 2 & 1.07 & - \\
& GREEDY POLYMAKE & - & - \\
\hline
\multirow{2}{0.2cm}{15} & COLGEN 2 & 0.01 & - \\
& GREEDY POLYMAKE & - & - \\
\hline
\multirow{2}{0.2cm}{16} & COLGEN 2 & $< 0.003$ & - \\
& GREEDY POLYMAKE & - & - \\
\hline
\end{tabular}
\caption{Average error value and volume of the approximate convex hull for $d = 8$.}
\label{tab:p2d8}
\end{table}

Regarding dimension 2 and 4, the column generation algorithm obtains lower average separation error in 5 cases over 6. Because the hyperplanes describing the convex hull have distance 0 from the set of positive points, we could expect an advantage of GREEDY POLYMAKE in terms of the volume of the resulting polytopes. Instead, COLGEN 2 obtains a better result in 2 cases over 6 and its volume values are close to those of GREEDY POLYMAKE in the remaining cases.

\section{Conclusions}
\label{sec:Conclusions}

We proposed a mathematical programming approach to a problem in computational geometry, which is to find a polyhedron with a limited number of faces, approximating the convex hull of a set of points (\CPS in the paper).

Ultimately, our approach allows us to design formulations for the \CPS, and ad-hoc algorithms which exploit column generation techniques, and rely on special purpose pricing routines.

The first task was to design formulations with a polynomial number of variables and constraints. Although our study began by inspirations from support vector machines, experiments indicate a direct adaptation not to be competitive. By elaborating on them, however, we could propose compact models which outperformed earlier MIPs from the literature. Key ideas were to properly change the penalty function, removing both the handling of margins and the hyperplane norm regularization terms.

The second task was to design more effective resolution algorithms. We used decomposition methods to obtain extended formulations and a column generation scheme. When solving pricing problems as MIPs, unfortunately, such a decomposition process had only a partial impact: even if smaller, pricing MIPs remain challenging to tackle with general purpose solvers. We designed an ad-hoc heuristic pricing routine that takes advantage from considering only a subset of points simultaneously and from parallelization. Its embedding in the column generation framework led to an effective algorithm for \CPS. In our experiments, in its best configuration, the column generation algorithm was able to solve instances up to dimension 8 with number of points up to over 10 thousand, comparing favourably to all the benchmark algorithms we have considered.

With respect to the target of obtaining convex hull approximations, our experiments highlight our approach to be effective. We mention three appealing features. First, our algorithms are flexible tools for numerical resolution; their structure makes it easy to tune them to be more accurate (at the price of higher computing times) or more quick (at the price of a few more separation errors). Standard approaches from the literature do not have such an option. Second, our models allow to possibly produce unbounded polyhedra as solutions; convex hull algorithms are forced to produce convex polytopes. Third, our models assume the number of hyperplanes, to be produced as faces, as part of data. This may actually become a hyperparameter in applications. That is, our methods allow parametric analyses to be performed, to evaluate the convex hull approximation error which is made by progressively allowing less faces, and therefore producing less complex polyhedra. It is interesting to note that, in our experiments, very few hyperplanes are enough to get approximations whose volume closely match that of the actual convex hull.

As a benchmark, we considered adaptations of convex hull algorithms from the literature, and in particular from the package Polymake; these latter show to be faster than our mathematical programming approaches when the number of dimensions is very low (providing approximations of similar quality), but get orders of magnitude slower as the dimensionality increases. The use of Polymake becomes inapplicable already in dimension $8$. We report that both approaches are affected by the number of points in the instance, which may however increase exponentially (in case the positive and negative points are produced by a proper sampling) or not (in case the positive and negative points are feasible and infeasible solutions of a MIP).

The set of  In our \CPS, the sets of positive and negative points are assumed to be given.
As perspective research, we plan to consider the case in which new positive and negative points can be generated by suitable queries, at a cost, as in \cite{Messana2024Actively}. We expect this line of research to yield methods which are applicable in even wider settings. \\

\noindent
{\bf Acknowledgments. }  The work has been partially supported by the University of Milan, Piano Sostegno alla Ricerca.

\bibliographystyle{informs2014}
\bibliography{references}

\begin{thebibliography}{33}
\providecommand{\natexlab}[1]{#1}
\providecommand{\url}[1]{\texttt{#1}}
\providecommand{\urlprefix}{URL }

\bibitem[{Amaldi et~al.(2012)Amaldi, Dhyani, \protect\BIBand{} Ceselli}]{amaldi12}
Amaldi E, Dhyani K, Ceselli A (2012) Column generation for the minimum hyperplanes clustering problem 25.

\bibitem[{Assarf et~al.(2017)Assarf, Gawrilow, Herr, Joswig, Lorenz, Paffenholz, \protect\BIBand{} Rehn}]{assarf2017computing}
Assarf B, Gawrilow E, Herr K, Joswig M, Lorenz B, Paffenholz A, Rehn T (2017) Computing convex hulls and counting integer points with polymake. \emph{Mathematical Programming Computation} 9:1--38.

\bibitem[{Astorino et~al.(2023)Astorino, Avolio, \protect\BIBand{} Fuduli}]{Astorino2023MaximummarginPS}
Astorino A, Avolio M, Fuduli A (2023) Maximum-margin polyhedral separation for binary multiple instance learning. \emph{EURO Journal on Computational Optimization} .

\bibitem[{Astorino et~al.(2021)Astorino, Francesco, Gaudioso, Gorgone, \protect\BIBand{} Manca}]{Astorino2021PolyhedralSV}
Astorino A, Francesco MD, Gaudioso M, Gorgone E, Manca B (2021) Polyhedral separation via difference of convex (dc) programming. \emph{Soft Computing} 25:12605 -- 12613.

\bibitem[{Astorino \protect\BIBand{} Fuduli(2015)}]{Astorino2015SupportVM}
Astorino A, Fuduli A (2015) Support vector machine polyhedral separability in semisupervised learning. \emph{Journal of Optimization Theory and Applications} 164:1039--1050.

\bibitem[{Astorino \protect\BIBand{} Gaudioso(2002)}]{Astorino2002PolyhedralST}
Astorino A, Gaudioso M (2002) Polyhedral separability through successive lp. \emph{Journal of Optimization Theory and Applications} 112:265--293.

\bibitem[{Balestriero et~al.(2022)Balestriero, Wang, \protect\BIBand{} Baraniuk}]{Balestriero2022DeepHullFC}
Balestriero R, Wang Z, Baraniuk R (2022) Deephull: fast convex hull approximation in high dimensions. \emph{ICASSP 2022 - 2022 IEEE International Conference on Acoustics, Speech and Signal Processing (ICASSP)} 3888--3892.

\bibitem[{Baum(1990)}]{baum1990polynomial}
Baum EB (1990) A polynomial time algorithm that learns two hidden unit nets. \emph{Neural Computation} 2(4):510--522.

\bibitem[{Baum(1991)}]{baum1991neural}
Baum EB (1991) Neural net algorithms that learn in polynomial time from examples and queries. \emph{IEEE Transactions on Neural Networks} 2(1):5--19.

\bibitem[{Blum et~al.(1995)Blum, Chalasani, Goldman, \protect\BIBand{} Slonim}]{blum1995learning}
Blum A, Chalasani P, Goldman SA, Slonim DK (1995) Learning with unreliable boundary queries. \emph{Proceedings of the eighth annual conference on Computational learning theory}, 98--107.

\bibitem[{{Bradford Barber} et~al.(1996){Bradford Barber}, Dobkin, \protect\BIBand{} Huhdanpaa}]{Barber1996TheQA}
{Bradford Barber} C, Dobkin DP, Huhdanpaa H (1996) The quickhull algorithm for convex hulls. \emph{ACM Trans. Math. Softw.} 22:469--483.

\bibitem[{Bykat(1978)}]{Bykat1978ConvexHO}
Bykat A (1978) Convex hull of a finite set of points in two dimensions. \emph{Inf. Process. Lett.} 7:296--298, \urlprefix\url{https://api.semanticscholar.org/CorpusID:31239431}.

\bibitem[{Carrizosa et~al.(2010)Carrizosa, Martin-Barragan, \protect\BIBand{} {Romero Morales}}]{carrizosa10}
Carrizosa E, Martin-Barragan B, {Romero Morales} D (2010) Binarized support vector machines. \emph{INFORMS Journal on Computing} 22(1):154--167.

\bibitem[{Eddy(1977)}]{Eddy1977ANC}
Eddy WF (1977) A new convex hull algorithm for planar sets. \emph{ACM Trans. Math. Softw.} 3:398--403, \urlprefix\url{https://api.semanticscholar.org/CorpusID:6464920}.

\bibitem[{Fajemisin et~al.(2024)Fajemisin, Maragno, \protect\BIBand{} {den Hertog}}]{fajemisin2024optimization}
Fajemisin AO, Maragno D, {den Hertog} D (2024) Optimization with constraint learning: A framework and survey. \emph{European Journal of Operational Research} 314(1):1--14.

\bibitem[{Gawrilow \protect\BIBand{} Joswig(2000)}]{gawrilow2000polymake}
Gawrilow E, Joswig M (2000) Polymake: a framework for analyzing convex polytopes. \emph{Polytopes—combinatorics and computation}, 43--73 (Springer).

\bibitem[{Gottlieb et~al.(2018)Gottlieb, Kaufman, Kontorovich, \protect\BIBand{} Nivasch}]{Gottlieb2018LearningCP}
Gottlieb LA, Kaufman E, Kontorovich A, Nivasch G (2018) Learning convex polytopes with margin. \emph{Neural Information Processing Systems}.

\bibitem[{Jayaram \protect\BIBand{} Fleyeh(2016)}]{jayaram2016convex}
Jayaram M, Fleyeh H (2016) Convex hulls in image processing: a scoping review. \emph{American Journal of Intelligent Systems} 6(2):48--58.

\bibitem[{Kantchelian et~al.(2014)Kantchelian, Tschantz, Huang, Bartlett, Joseph, \protect\BIBand{} Tygar}]{Kantchelian2014LargeMarginCP}
Kantchelian A, Tschantz MC, Huang L, Bartlett PL, Joseph AD, Tygar JD (2014) Large-margin convex polytope machine. \emph{Neural Information Processing Systems}.

\bibitem[{Klivans et~al.(2004)Klivans, O'Donnell, \protect\BIBand{} Servedio}]{klivans2004learning}
Klivans AR, O'Donnell R, Servedio RA (2004) Learning intersections and thresholds of halfspaces. \emph{Journal of Computer and System Sciences} 68(4):808--840.

\bibitem[{Kwek \protect\BIBand{} Pitt(1996)}]{kwek1996pac}
Kwek S, Pitt L (1996) Pac learning intersections of halfspaces with membership queries. \emph{Proceedings of the ninth annual conference on Computational learning theory}, 244--254.

\bibitem[{Leng et~al.(2020)Leng, He, Liu, ping Qin, \protect\BIBand{} Li}]{Leng2020ASC}
Leng Q, He Z, Liu Y, ping Qin Y, Li Y (2020) A soft-margin convex polyhedron classifier for nonlinear task with noise tolerance. \emph{Applied Intelligence} 51:453--466.

\bibitem[{Lombardi \protect\BIBand{} Milano(2018)}]{lombardi2018boosting}
Lombardi M, Milano M (2018) Boosting combinatorial problem modeling with machine learning. \emph{Proceedings of the Twenty-Seventh International Joint Conference on Artificial Intelligence} .

\bibitem[{Manwani \protect\BIBand{} Sastry(2010)}]{Manwani2010LearningPC}
Manwani N, Sastry PS (2010) Learning polyhedral classifiers using logistic function. \emph{Asian Conference on Machine Learning}.

\bibitem[{Megiddo(1988)}]{megiddo88}
Megiddo N (1988) On the complexity of polyhedral separability. \emph{Discrete and Computational Geometry} 3.

\bibitem[{Messana(2024)}]{RD_UNIMI/IE2MX0_2024}
Messana R (2024) {Replication Data for: Mathematical programming algorithms for convex hull approximation with a hyperplane budget}. \urlprefix\url{http://dx.doi.org/10.13130/RD_UNIMI/IE2MX0}.

\bibitem[{Messana et~al.(2024)Messana, Chen, \protect\BIBand{} Lodi}]{Messana2024Actively}
Messana R, Chen R, Lodi A (2024) Actively learning combinatorial optimization using a membership oracle. \emph{ArXiv} abs/2405.14090.

\bibitem[{Orsenigo \protect\BIBand{} Vercellis(2007)}]{Orsenigo2007AccuratelyLF}
Orsenigo C, Vercellis C (2007) Accurately learning from few examples with a polyhedral classifier. \emph{Computational Optimization and Applications} 38:235--247.

\bibitem[{Pawlak \protect\BIBand{} Krawiec(2017)}]{Pawlak2017AutomaticSO}
Pawlak TP, Krawiec K (2017) Automatic synthesis of constraints from examples using mixed integer linear programming. \emph{Eur. J. Oper. Res.} 261:1141--1157.

\bibitem[{Raviv et~al.(2018)Raviv, Hazan, \protect\BIBand{} Osadchy}]{Raviv2018HingeMinimaxLF}
Raviv D, Hazan T, Osadchy M (2018) Hinge-minimax learner for the ensemble of hyperplanes. \emph{J. Mach. Learn. Res.} 19:62:1--62:30.

\bibitem[{Schede et~al.(2019)Schede, Kolb, \protect\BIBand{} Teso}]{Schede2019LearningLP}
Schede E, Kolb S, Teso S (2019) Learning linear programs from data. \emph{2019 IEEE 31st International Conference on Tools with Artificial Intelligence (ICTAI)} 1019--1026.

\bibitem[{Strekalovsky et~al.(2015)Strekalovsky, Gruzdeva, \protect\BIBand{} Orlov}]{Strekalovsky2015OnTP}
Strekalovsky AS, Gruzdeva TV, Orlov AV (2015) On the problem polyhedral separability: a numerical solution. \emph{Automation and Remote Control} 76:1803--1816.

\bibitem[{Vempala(1997)}]{vempala1997random}
Vempala S (1997) A random sampling based algorithm for learning the intersection of half-spaces. \emph{Proceedings 38th Annual Symposium on Foundations of Computer Science}, 508--513 (IEEE).

\end{thebibliography}

\newpage

\APPENDIX{}

\section{Classical SVM Models and Adaptations}
\label{sect:svm-model}

In order to introduce our MIPs, we report below a classical soft margin SVM model shown below.

\begin{align}
\text{min } & \sum_{i \in I} C_i \xi_i - \sum_{j \in J} w_j^2 \label{obj:svm} \\
\text{s.t. } & b + \sum_{j \in J} w_j a_{ij} \geq 1 - \xi_i & \forall i \in I^+ \label{st:svm_positive} \\
            & b + \sum_{j \in J} w_j a_{ij} \leq - 1 + \xi_i & \forall i \in I \setminus I^+ \label{st:svm_negative} \\
            & \xi_i \geq 0 & \forall i \in I \label{st:svm_domain_xi} \\
            & b \in \mathbb{R} \label{st:svm_domain_b} \\
            & w_j \in \mathbb{R} & \forall j \in J \label{st:svm_domain_w}
\end{align}

Variables $\boldsymbol{w}$ encode hyperplane coefficients; variables $\xi_i$ encode misclassification error on point $i \in I$.
Constraints \eqref{st:svm_positive} and \eqref{st:svm_negative} impose correct classification of the points in both positive and negative classes, with the possibility of misclassifying some of the points paying in terms of the corresponding violation variables. The objective is to minimize the sum of the violations minus a regularization term proportional to the squared 2-norm of the hyperplane coefficients $\boldsymbol{w}$. Each constants $C_i \geq 0$ is a weight that quantifies the importance of correctly classifing the $i^{th}$ point. In practice, all such weights are often set equal to a same constant $C$.\\
We designed the following adaptation to \CPS.
\begin{align}
\text{min } & \sum_{i \in I^-} \hat c(\boldsymbol{a}^-_i) \xi_i^k x_i^k \label{obj:model_a} \\
\text{s.t. } & \sum_{k \in K} x_i^k = 1 & \forall i \in I^- \label{st:model_a_x} \\
            & b^k + \sum_{j \in J} w_j^k a_{ij} \geq 1 & \forall i \in I^+, \forall k \in K \label{st:model_a_positive} \\
            & b^k + \sum_{j \in J} w_j^k a_{ij} \leq - 1 + \xi_i^k & \forall i \in I^-, \forall k \in K \label{st_model_a_negative} \\
            & \xi_i^k \geq 0 & \forall i \in I^-, \forall k \in K \label{st:model_a_domain_xi} \\
            & x_i^k \in \{0, 1\} & \forall i \in I^-, \forall k \in K \label{st:model_a_domain_x} \\
            & b^k \in \mathbb{R} & \forall k \in K \label{st:model_a_domain_b} \\
            & w_j^k \in \mathbb{R} & \forall j \in J, \forall k \in K \label{st:model_a_domain_w}
\end{align}
The model prescribes to find $\mathcal{K}$ hyperplane equations, encoded by variables $w_j^k$ with $j \in J$ and $k \in K$. A set of non-negative variables $\xi_i^k$ encode the distance from point $i \in I^-$ to hyperplane $k$, if $i$ is on the wrong side of the hyperplane $k$; they take value $0$ otherwise, according to constraints \eqref{st_model_a_negative}.
Since \CPS requires every hyperplane in a solution to be valid, we neglect the violation variables corresponding to positive points, imposing them to be always on the correct side of each hyperplane, according to constraints \eqref{st:model_a_positive}. 
We introduce binary variables $x_i^k$ for $i \in I^-$ and $k \in K$. When $x_i^k = 1$, it means that point $i$ is assigned to the $k^{th}$ hyperplane. According to constraints \eqref{st:model_a_x}, each negative point must be assigned to exactly one hyperplane. Constraints \eqref{st:model_a_positive} force the hyperplanes to be valid and constraints \eqref{st_model_a_negative} require the hyperplanes to keep the negative points on a same side with possibility of violation. The objective function sums up the weighted violations over the negative points with respect to the hyperplanes they are assigned to, where the weights are expressed as values of a function $\hat c: \mathbb{R}^d \to [0, +\infty)$. In force of constraints \eqref{st:model_a_x}, if we fix the $\mathcal{K}$ hyperplanes, the minimum value of the objective function is the sum over the negative points of the minimum weighted violation with respect to the hyperplanes.

\section{Column Generation and Hyperplane Choice}
\label{sec:ColumnGenerationAndHyperplaneChoice}

\begin{algorithm}
\small
\caption{Column Generation and Hyperplane Choice}\label{alg:CGandHC}
\begin{algorithmic}
    \State \textbf{Input:} point limit $thr$ and maximum number of processes $N^{max}$ (only when using AHP, see Section \ref{subsec:ahp})
    \State
    \State $\bar{N} \gets 0$
    \State $R \gets -\infty$
    \State $E \gets +\infty$
    \State $\boldsymbol{x}^{prev} \gets \boldsymbol{1} \in \mathbb{N}^n$
    \State Set $\boldsymbol{W}$ to empty list
    \State Set $\boldsymbol{b}$ to empty list
    \State Initialize RMM without $z$ variables
    \While{$R < 0$ and $E > 0$}
        \State Solve RMM and get dual variable values $\boldsymbol{\lambda} = (\lambda_i)_{i \in I^-}$ and $\mu$
        \State Solve PM($\boldsymbol{\lambda}$, $\mu$) (or AHP($\boldsymbol{\lambda}$, $\mu$, $thr$, $\boldsymbol{x}^{prev}$, $N^{max}$), see Section \ref{subsec:ahp}) and get:
        \State\ \ \ (a) the number $N > 0$ of columns generated
        \State\ \ \ (b) vectors $\boldsymbol{x}^k \in \{0, 1\}^n$ for each $k \in \{1, \ldots, N\}$
        \State\ \ \ (c) scalars $b^k \in \mathbb{R}$ for each $k \in \{1, \ldots, N\}$
        \State\ \ \ (d) vectors $\boldsymbol{w}^k \in \mathbb{R}^d$ for each $k \in \{1, \ldots, N\}$
        \State $\boldsymbol{r} \gets \left(- \sum_{i \in I^-} x_i^k \cdot \lambda_i + \mu\right)_{k=1}^{N}$
        \State $k^{min} \gets argmin_{k=1}^{N} r_k$
        \State $R \gets r_{k^{min}}$
        \State $\boldsymbol{x}^{prev} \gets \boldsymbol{x}^{k^{min}}$
        \State Add variable $z^{\bar{N} + k}$ and column $\boldsymbol{c}^{\bar{N} + k} = \left(x^k \mathbin\Vert (-1)\right)^{\top}$ to RMM for each $k \in \{1, \ldots, N\}$ such that $r_k < 0$
        \State Solve HCP and get the objective value, i.e., the error value $E$
        \State Append $b^k$ to $\boldsymbol{b}$ for every $k \in \{1, \ldots, N\}$
        \State Append $\boldsymbol{w}^k$ to $\boldsymbol{W}$ for every $k \in \{1, \ldots, N\}$
        \State $\bar{N} \gets \bar{N} + N$
    \EndWhile \\
    \Return $\{(\boldsymbol{b}(k),\ \boldsymbol{W}(k))\ |\ k \in \{1, \ldots, \bar{N}\} \ \text{and}\ z^k = 1\ \text{in the last HCP solution}\}$
\end{algorithmic}
\end{algorithm}

\section{Adaptation to \CPS of a MIP from Literature}
\label{sec:AppendixOV2007}

We show below the adaptation for \CPS of the model proposed in \citep*{Orsenigo2007AccuratelyLF} for convex polyhedral separation.  We neglected part of its objective function together with the respective variables and the related constraints, in order to make the function equal to the separation error of \CPS. We also added constraints \eqref{constr:ov_d} to impose correct separation of all the positive points. The parameter $\epsilon$ is greater than $0$, while $Q$ is equal to $|J|D$, where $D$ is greater or equal than the diameter of $\mathcal{D}^+ \cup \mathcal{D}^-$.

The binary variables take value 1 according to the following conditions:
\begin{enumerate}
    \item $e_i$ if and only if point $i$ is misclassified with respect to its true label;
    \item $d_i^k$ if and only if point $i$ satisfies $\boldsymbol{w}^k \cdot \boldsymbol{a}_i \geq w_0^k$, i.e. hyperplane $k$ classifies point $i$ as ``positive'';
    \item $m_i$ if and only if point $i$ lies outside the convex region identified by the $\mathcal{K}$ hyperplanes, i.e. $d_i^k = 1$ for some $k \in K$;
    \item $v$ if and only if the points inside the polyhedron identified by the $\mathcal{K}$ hyperplanes are labeled as ``1'' by the model;
\end{enumerate}
The continuous variables represent the parameters of the hyperplanes. The objective function of the model is a weighted sum of the separation errors over all the points, positive and negative. Constraints \eqref{constr:ov_classification_1} and \eqref{constr:ov_classification_2} identify, for each point $i$, which halfspace determined by hyperplane $k$ the point is in. Constraints \eqref{constr:ov_m_upper} ensures that $m_i = 0$ if $d_i^k = 0$ for every $k \in K$, while constraints \eqref{constr:ov_m_lower} impose that $m_i = 1$ if $d_i^k = 1$ for some $k \in K$. Constraints \eqref{constr:ov_m_e_1} and \eqref{constr:ov_m_e_2} characterize the error values $e_i$ in terms of $m_i$, $v$ and the true labels of the points. Finally, constraints \eqref{constr:ov_d} are not present in the original model \citep*{Orsenigo2007AccuratelyLF} and impose that all the points whose true label is 1 are classified as positive by all the hyperplanes. These constraints make the model suitable to solve PCAB, which requires all the positive samples to be classified correctly.
\begin{align}
\text{min } & \sum_{i \in I} c(\boldsymbol{a}_i) e_i \\
\text{s.t. } & b^k + \sum_{j \in J} w_j^k a_{ij} \geq - Q d_i^k & \forall i \in I, \forall k \in K \label{constr:ov_classification_1} \\
            & b^k + \sum_{j \in J} w_j^k a_{ij} \leq \left( 1 - d_i^k \right) Q - \epsilon & \forall i \in I, \forall k \in K \label{constr:ov_classification_2} \\
            & m_i \leq \sum_{k \in K} d_i^k & \forall i \in I \label{constr:ov_m_upper} \\
            & |K| m_i \geq \sum_{k \in K} d_i^k & \forall i \in I \label{constr:ov_m_lower} \\
            & 2 v - 1 - \ell(\boldsymbol{a}_i) (1 - 2 m_i) \leq 2 e_i & \forall i \in I \label{constr:ov_m_e_1} \\
            & 2 v - 1 - \ell(\boldsymbol{a}_i) (1 - 2 m_i) \geq - 2 e_i & \forall i \in I \label{constr:ov_m_e_2} \\
            & d_i^k = 0 & \forall i \in I^+, \forall k \in K \label{constr:ov_d} \\
            & e_i \in \{0, 1\} & \forall i \in I \\
            & d_i^k \in \{0, 1\} & \forall i \in I, \forall k \in K \\
            & m_i \in \{0, 1\} & \forall i \in I \\
            & v \in \{0, 1\} \\
            & b^k \in \mathbb{R} & \forall k \in K \\
            & w_j^k \in \mathbb{R} & \forall j \in J, \forall k \in K
\end{align}

\section{Data Generation}
\label{sec:DataGeneration}

\begin{table}[t]
\small
\centering
\begin{tabular}{ c r r|r r|r r }
 \thead{$d$} & \thead{N. det.\\ positive points} & \thead{N. det.\\ negative points} & \thead{N. random \\ positive points} & \thead{N. random \\ negative points} & \thead{Total n. \\ positive points} & \thead{Total n. \\ negative points} \\
 \hline
 2 & 4 & 8 & 141 & 200 & 145 & 208 \\
 4 & 16 & 64 & 200 & 500 & 216 & 564 \\
 8 & 256 & 2048 & 282 & 8000 & 538 & 10048 \\
\end{tabular}
\caption{Details of the datasets for each $d \in \{2, 4, 8\}$.}
\label{tab:hypercubicData}
\end{table}

For the first part of the experiments, we designed and generated what we call ``hypercubic data'', i.e. a dataset intended with the goal of obtaining instances whose minimum separation budget is known \emph{a priori}. We refer to this dataset as D1 in the paper. The polyhedron to be approximated is the $d$-dimensional hypercube $H_d = [0, 1]^d$. After fixing a border gap $\gamma = 0.04$, for each vertex $v$ of the hypercube, we place:

\begin{enumerate}
    \item a positive point obtained from $v$ by replacing its 0 components with $\gamma$ and its 1 components with $1 - \gamma$;
    \item $d$ negative points, obtained from $v$ by applying two subsequent replacements. With the first one, we replace every 0 with $- \gamma$ and every 1 with $1 + \gamma$. Then, in turn for each $j \in \{1, 2, \ldots, d\}$, we replace the $j^{th}$ element. If it is $- \gamma$, it becomes $2 \gamma$ and if it is $1 + \gamma$ it becomes $1 - 2 \gamma$.
\end{enumerate}

\begin{figure}[t]
     \centering
     \begin{subfigure}[b]{0.49\textwidth}
         \centering
         \includegraphics[width=\textwidth]{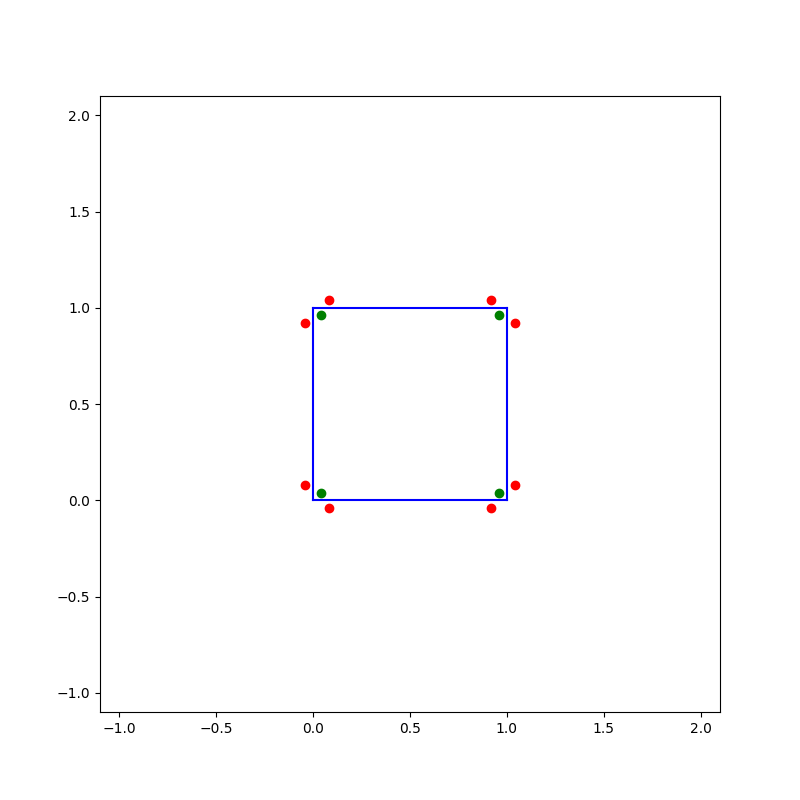}
         \caption{Corner points}
         \label{fig:cornerPoints}
     \end{subfigure}
     \hfill
     \begin{subfigure}[b]{0.49\textwidth}
         \centering
         \includegraphics[width=\textwidth]{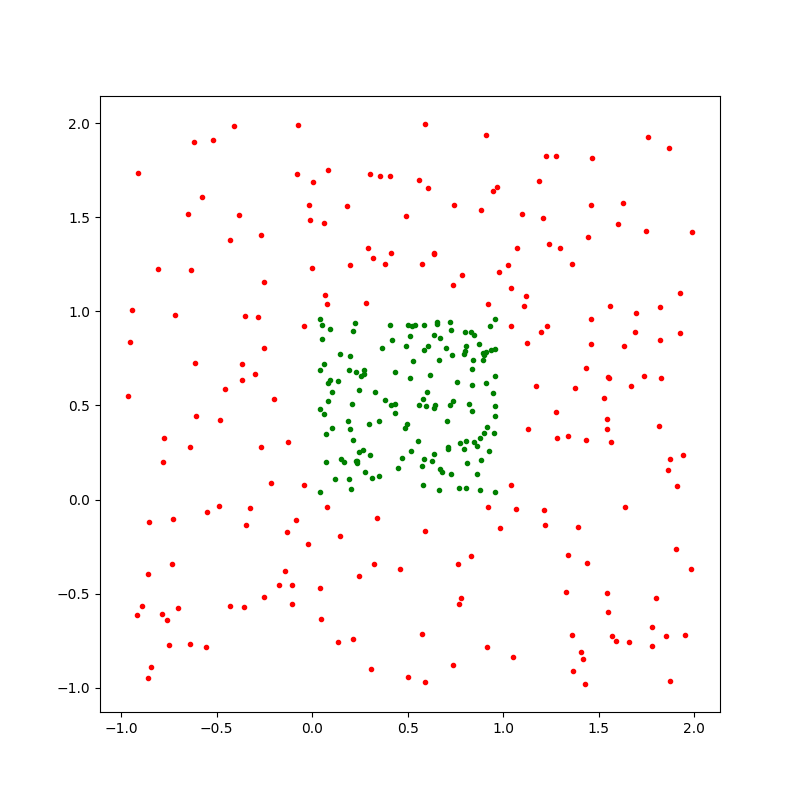}
         \caption{Example of instance}
         \label{fig:hypercubicDataInstance}
     \end{subfigure}
        \caption{Hypercubic data in dimension 2: points placed deterministically (\ref{fig:cornerPoints}) and entire instance representation including random points (\ref{fig:hypercubicDataInstance}).}
        \label{fig:hypercubicData}
\end{figure}

The result of this generation for $d = 2$ is shown in Fig.~\ref{fig:cornerPoints}. The points are placed in such a way that the minimum separation budget is $2d$. In order to enlarge this initial set, random points are added, both positive and negative (Fig.~\ref{fig:hypercubicDataInstance}), with the negative ones lying inside $[-1, 2]^d \setminus H^d$. We generated data in dimension 2, 4 and 8, producing 100 instances with the same number of points of each type, for each of the three dimensions. The number of deterministic, random and total points is shown in Table \ref{tab:hypercubicData}.

The dataset for the second part of the experiments, that we call D2, is generated as D1, but without any border gap ($\gamma = 0$) and without imposing any corner points in the instances. In this way, the minimum separating budget is possibly less that $2d$, while the number of facets of the convex hull of the positive points can be greater than $2d$. The number of positive and negative points in each dataset is the same as in the counterpart with corner points excluding such points.

\newpage

\section{Additional Figures and Tables}
\label{sec:AdditionalFiguresAndTables}

\begin{figure}[h]
    \centering
    \includegraphics[width=0.63\textwidth]{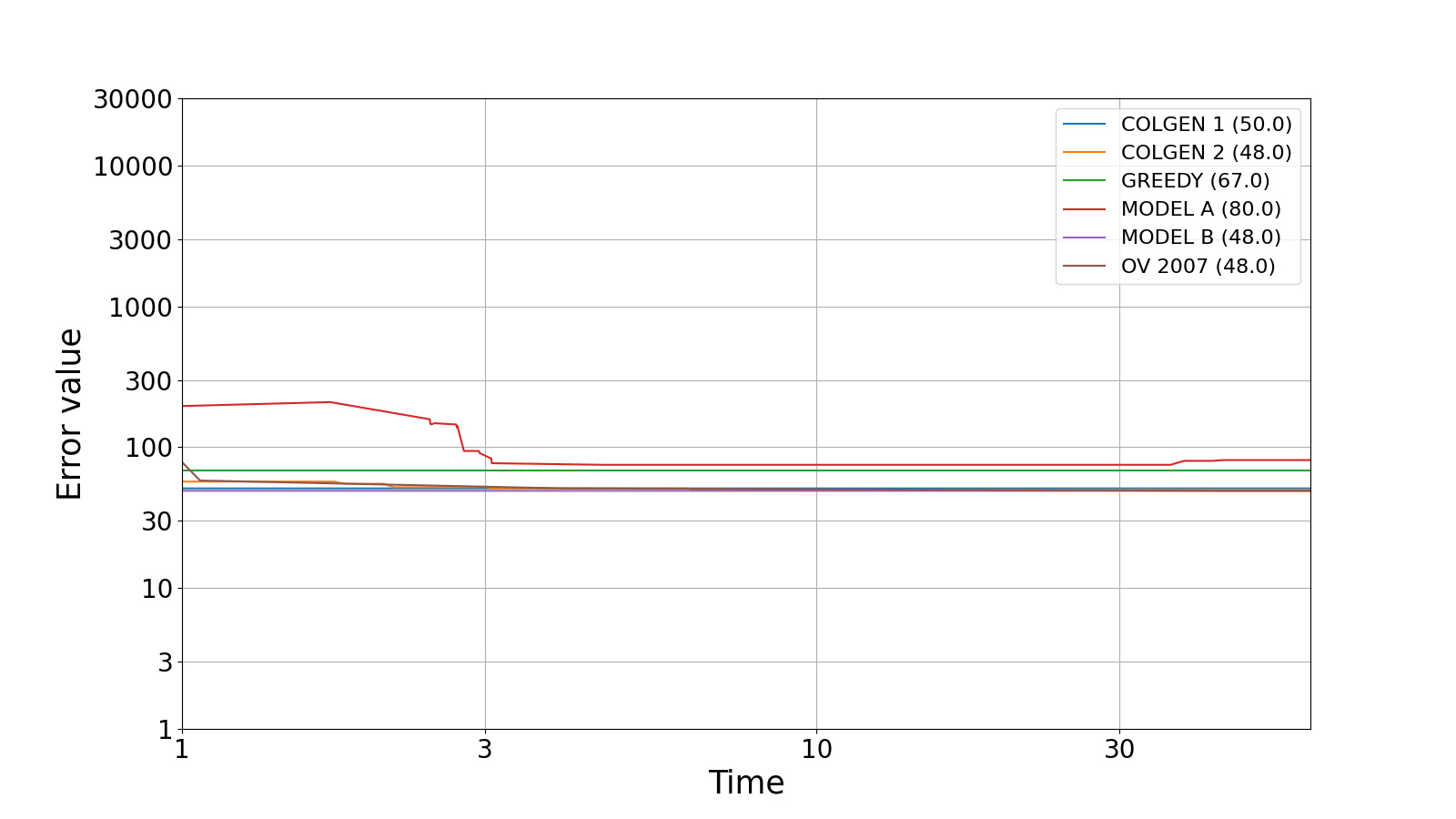}
    \caption{Error value with 2-dimensional hypercubic data: $\mathcal{K} = 2$.}
    \label{fig:d2K2}
\end{figure}

\begin{figure}[h]
    \centering
    \includegraphics[width=0.63\textwidth]{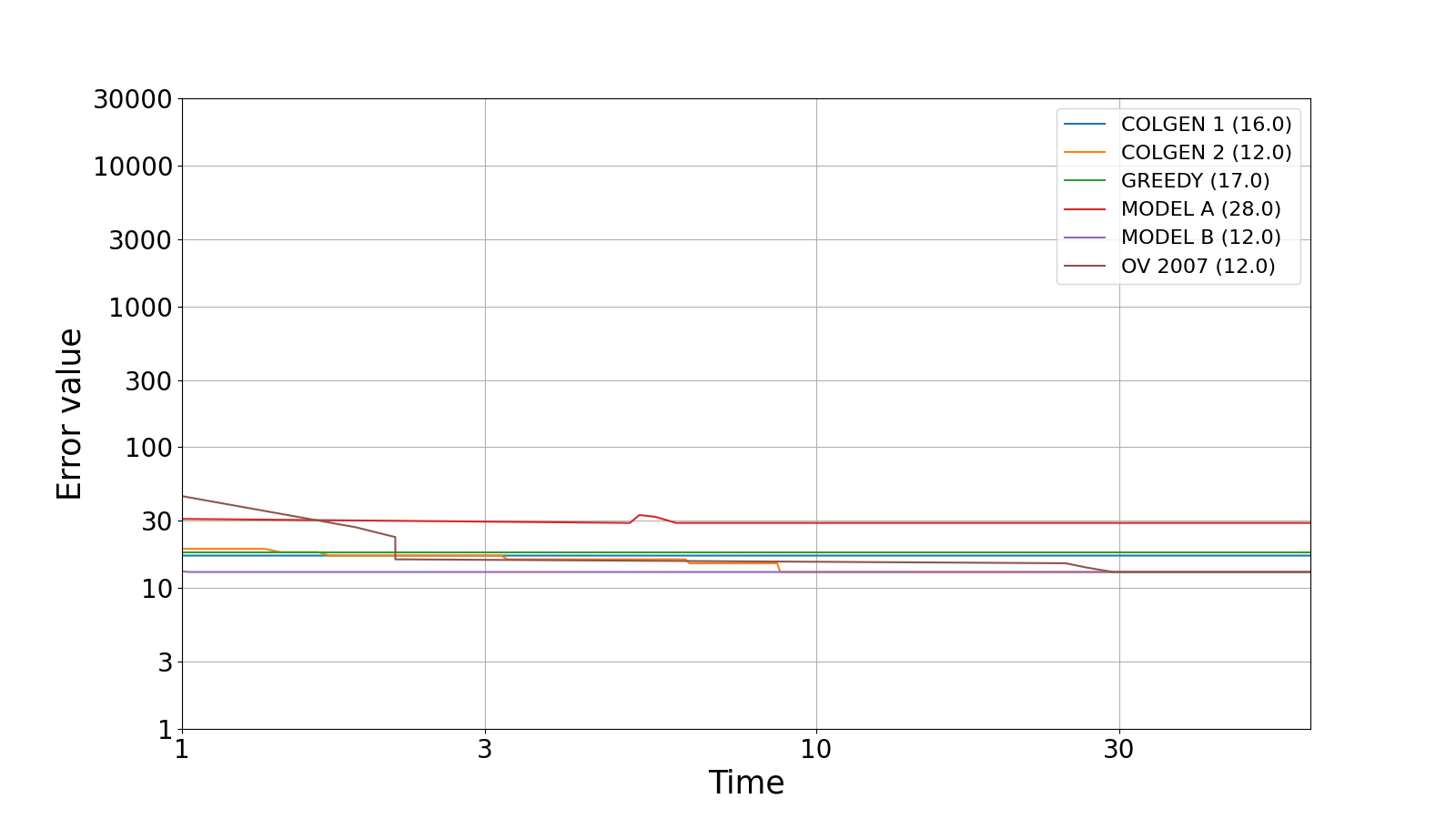}
    \caption{Error value with 2-dimensional hypercubic data: $\mathcal{K} = 3$.}
    \label{fig:d2K3}
\end{figure}

\begin{figure}[h]
    \centering
    \includegraphics[width=0.63\textwidth]{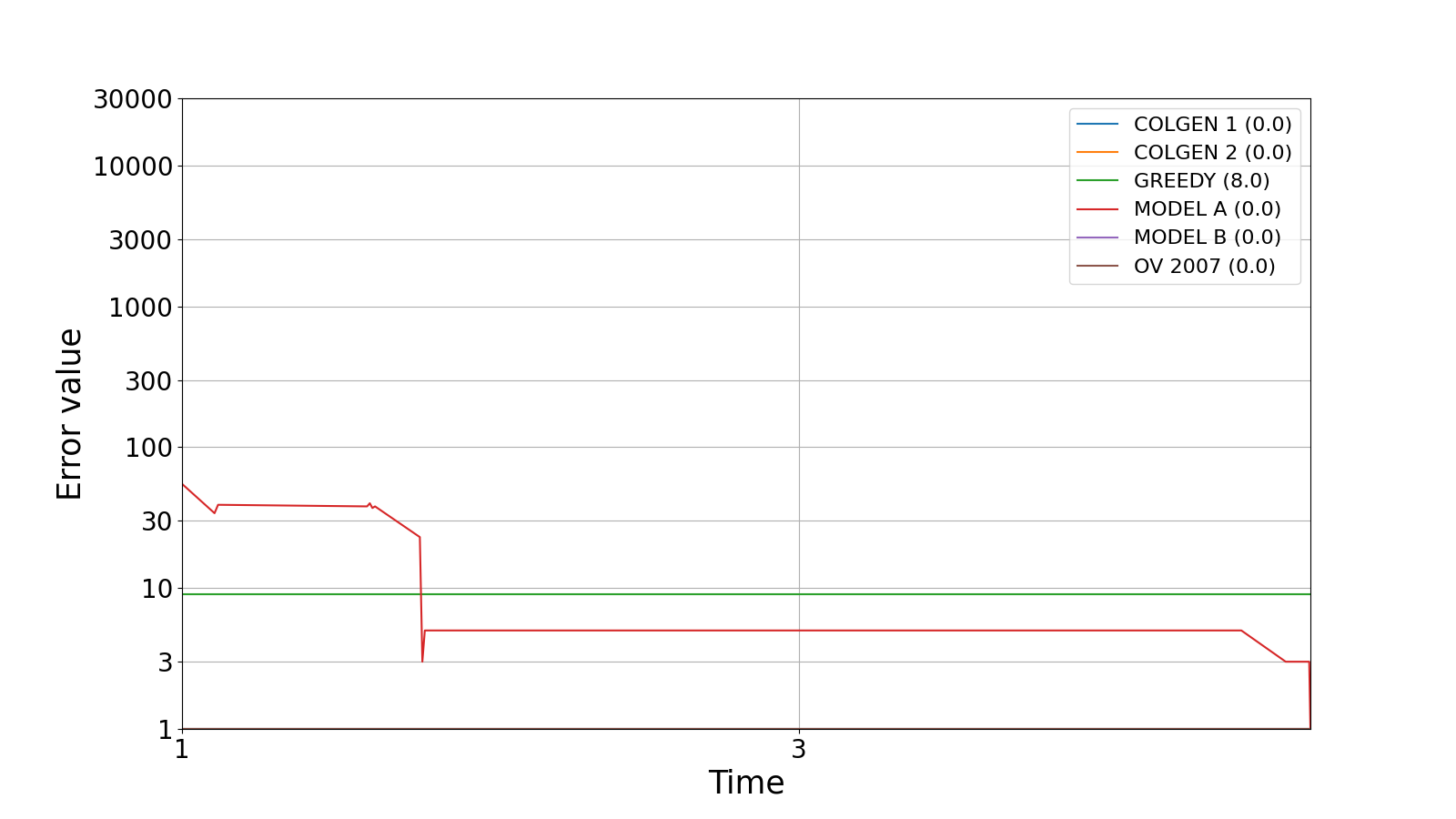}
    \caption{Error value with 2-dimensional hypercubic data: $\mathcal{K} = 4$.}
    \label{fig:d2K4}
\end{figure}

\begin{figure}[h]
    \centering
    \includegraphics[width=0.63\textwidth]{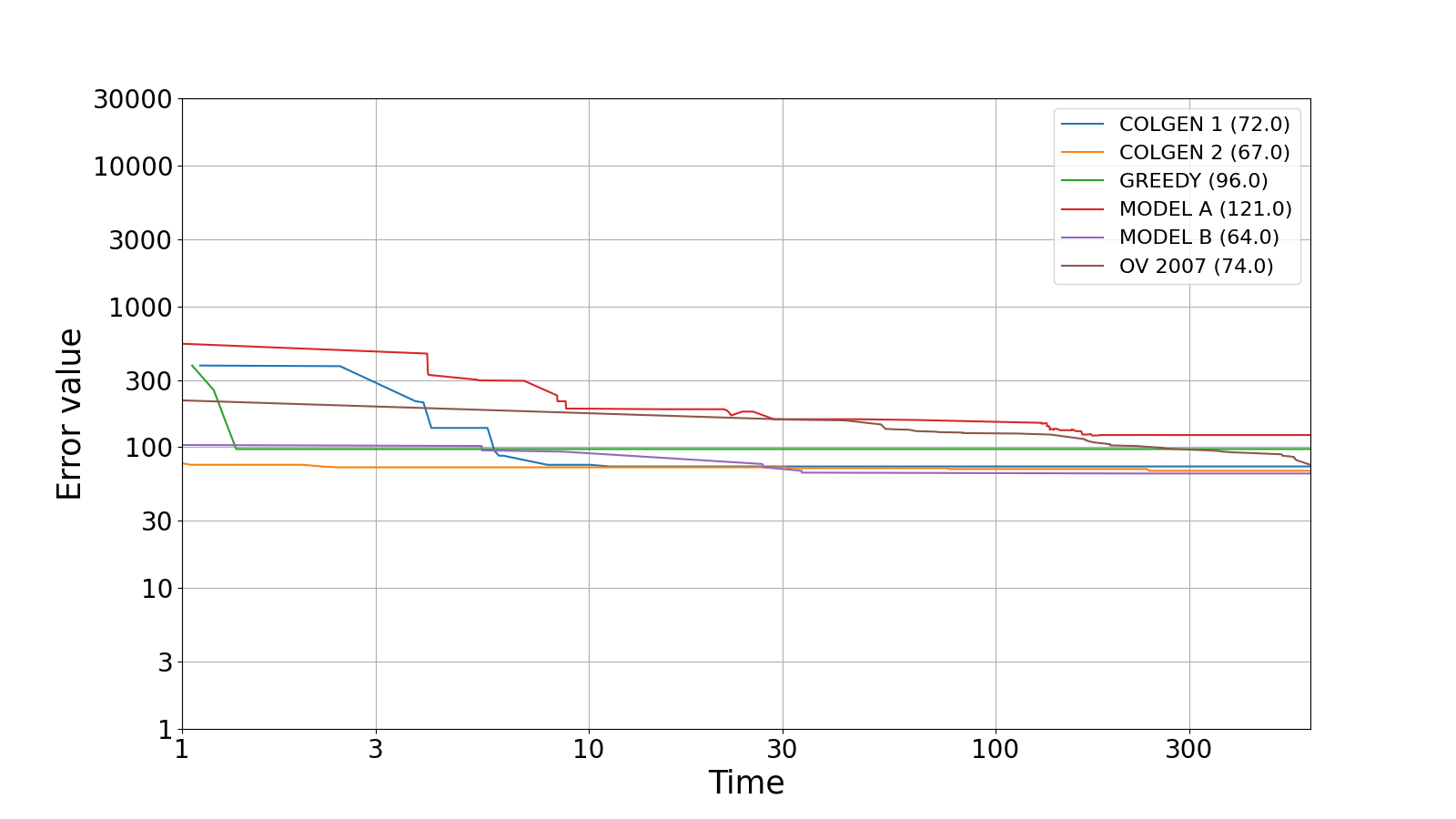}
    \caption{Error value with 4-dimensional hypercubic data: $\mathcal{K} = 4$.}
    \label{fig:d4K4}
\end{figure}

\begin{figure}[h]
    \centering
    \includegraphics[width=0.63\textwidth]{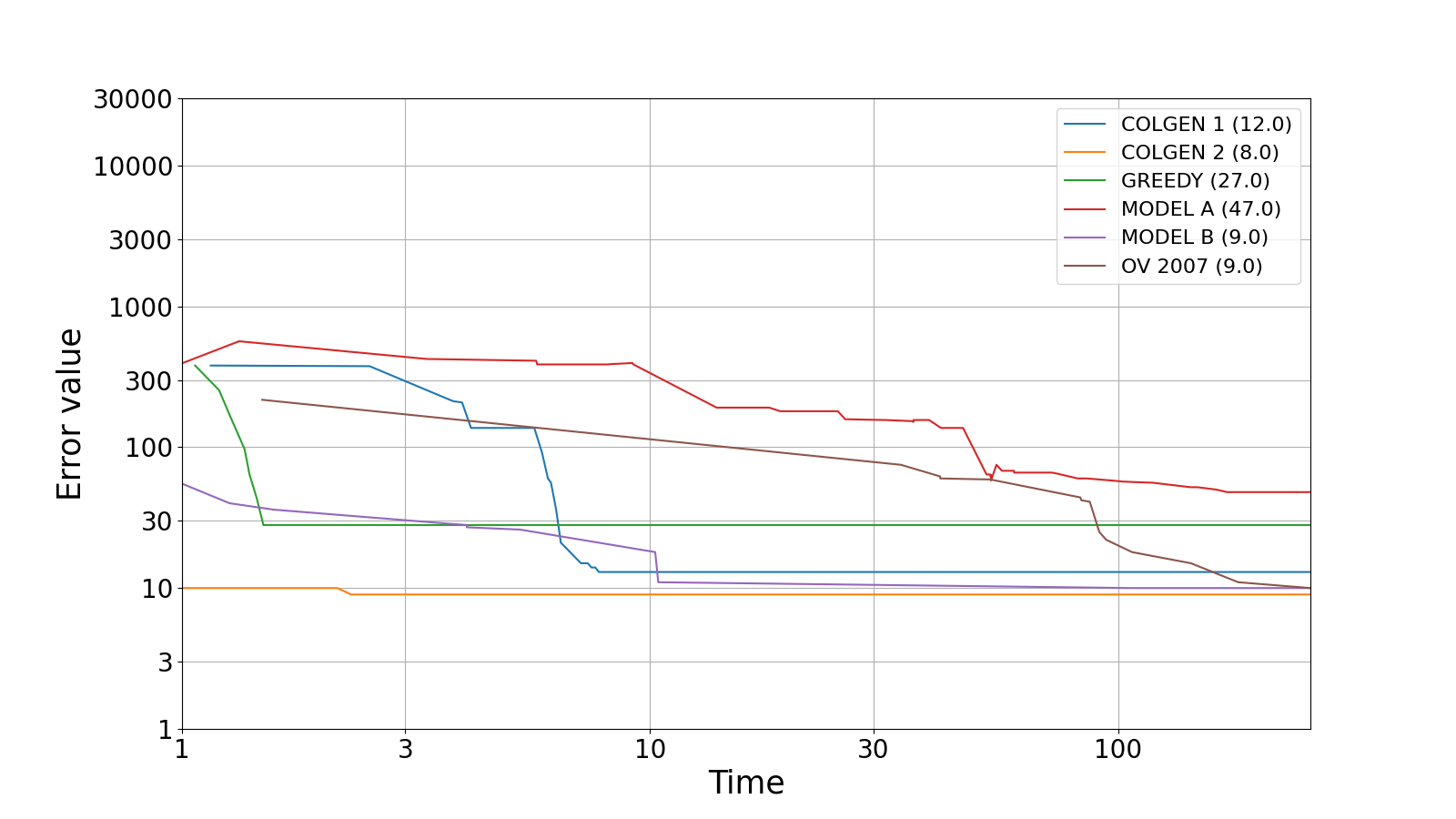}
    \caption{Error value with 4-dimensional hypercubic data: $\mathcal{K} = 7$.}
    \label{fig:d4K7}
\end{figure}

\begin{figure}[h]
    \centering
    \includegraphics[width=0.63\textwidth]{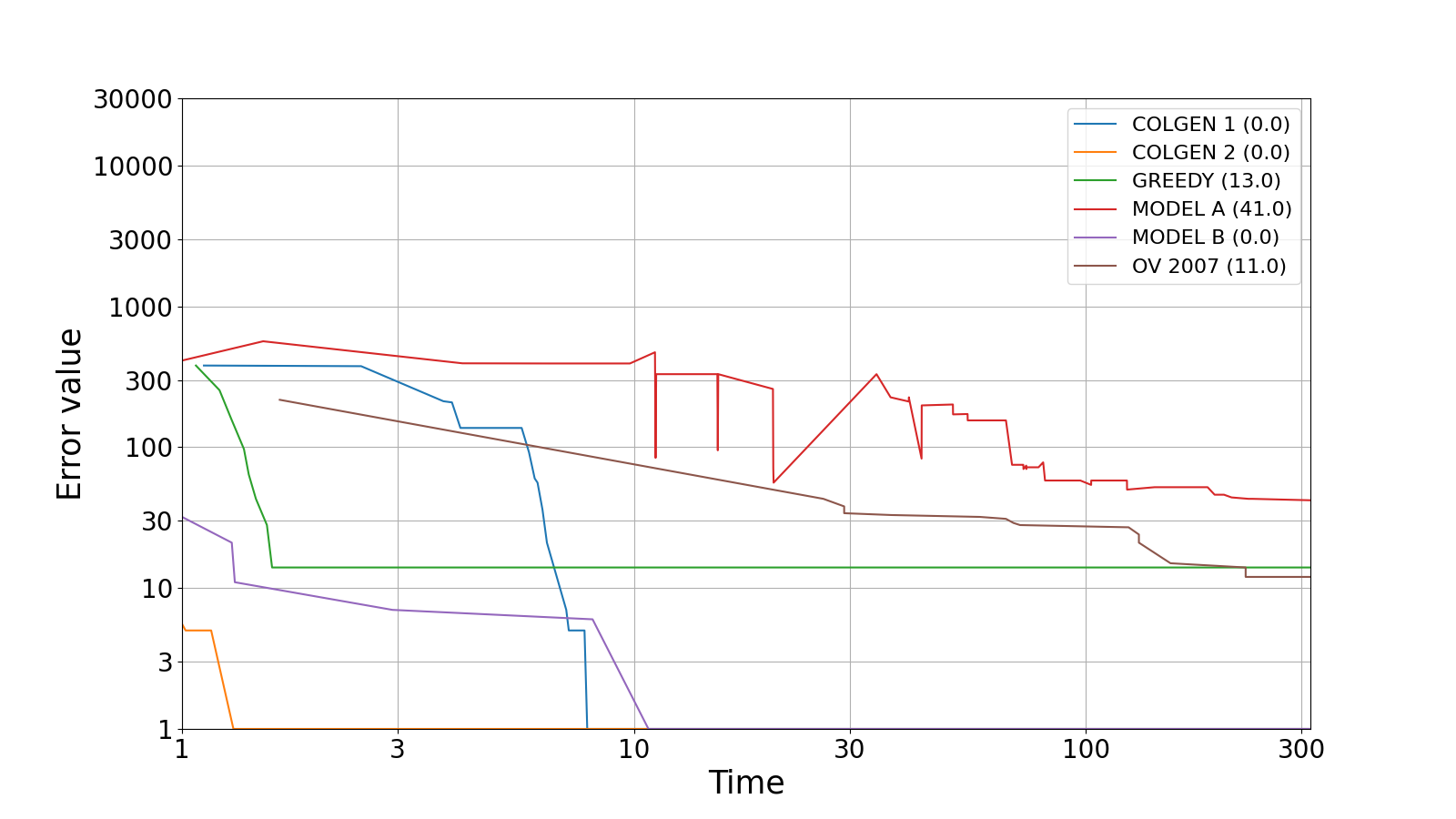}
    \caption{Error value with 4-dimensional hypercubic data: $\mathcal{K} = 8$.}
    \label{fig:d4K8}
\end{figure}

Figures \ref{fig:d2K2} - \ref{fig:d4K8} present the results of preliminary experiments on representative instances with corner points, both in dimension 2 (Figures \ref{fig:d2K2} - \ref{fig:d2K4}) and in dimension 4 (Figures \ref{fig:d4K4} - \ref{fig:d4K8}). Such Figures are analogous to Figures \ref{fig:d8K8} - \ref{fig:d8K16} presented in Section \ref{sec:ExperimentalSettingAndResults} for dimension 8. In dimension 2, the performance difference between different methods is scarse. In dimension 4, COLGEN 1, COLGEN 2 and MODEL B confirm to obtain the best results in terms of final error value. \\

Tables \ref{tab:massive_colgen_2_1} - \ref{tab:massive_colgen_8_1} show numerical results specific to our column generation algorithms, COLGEN 1 and COLGEN 2. We report number of iterations, number of columns generated and last plain pricing gap between objective and lower bound (only for for COLGEN 1). We also show the amount of time each of the two algorithms spent solving the pricing problem, the restricted master problem and the hyperplane choice problem respectively. We observe that the number of iterations that COLGEN 2 is able to complete can be significantly lower than that reported by COLGEN 1 when $\mathcal{K} < 2d$, especially with the datasets in higher dimensions.

However, the parallelism inherent in COLGEN 2 enables generation of a number of columns that is often greater with respect to the number of COLGEN 1. The number of parallel processes is an essential hyperparameter in this sense: extra preliminary experiments with a single process have shown far worse results in terms of error values compared to our standard experiments. The values of the last pricing gap of COLGEN 1 decrease when $\mathcal{K}$ increases, which indicates that finding a solution to the pricing model is especially harder for the solver when the hyperplane budget is lower. Finally, we can notice that the majority of the running times of the column generation algorithms is spent in solving the pricing model, while the restricted master model takes the smallest fraction of time and the hyperplane choice problem is in between.

\newpage

\newpage

\begin{table}[h!]
\small
\centering
\begin{tabular}{c l r r r r r r}
$\mathcal{K}$ & method & \makecell{n. iterations} & \makecell{n. columns} & \makecell{last pricing gap} & \makecell[r]{PM time\\(seconds)} & \makecell[r]{RMM time\\(seconds)} & \makecell[r]{HCM time\\(seconds)} \\
\hline
\multirow{2}{0.2cm}{2} & COLGEN 1 & 71.69 & 71.69 & 74.53 & 5.797 & 0.078 & 0.130 \\
& COLGEN 2 & 65.09 & 520.72 & - & 5.289 & 0.165 & 0.574 \\
\hline
\multirow{2}{0.2cm}{3} & COLGEN 1 & 77.09 & 77.09 & 34.21 & 5.759 & 0.097 & 0.145 \\
& COLGEN 2 & 65.07 & 520.55 & - & 5.271 & 0.169 & 0.590 \\
\hline
\multirow{2}{0.2cm}{4} & COLGEN 1 & 33.91 & 33.91 & 0.46 & 1.242 & 0.057 & 0.073 \\
& COLGEN 2 & 5.35 & 42.79 & - & 0.437 & 0.012 & 0.019 \\
\hline
\end{tabular}
\caption{Column generation results for $d = 2$: number of iterations, number of columns generated, last pricing gap and running times (pricing model, restricted master model, hyperplane choice model).}
\label{tab:massive_colgen_2_1}
\end{table}

\begin{table}[h!]
\small
\centering
\begin{tabular}{c l r r r r r r}
$\mathcal{K}$ & method & \makecell{n. iterations} & \makecell{n. columns} & \makecell{last pricing gap} & \makecell[r]{PM time\\(seconds)} & \makecell[r]{RMM time\\(seconds)} & \makecell[r]{HCM time\\(seconds)} \\
\hline
\multirow{2}{0.2cm}{4} & COLGEN 1 & 334.37 & 334.37 & 114.322 & 53.185 & 1.099 & 5.633 \\
& COLGEN 2 & 110.55 & 884.40 & - & 15.109 & 0.818 & 44.449 \\
\hline
\multirow{2}{0.2cm}{7} & COLGEN 1 & 391.23 & 391.23 & 54.84 & 53.072 & 1.269 & 4.445 \\
& COLGEN 2 & 280.72 & 2245.76 & - & 36.563 & 2.780 & 20.506 \\
\hline
\multirow{2}{0.2cm}{8} & COLGEN 1 & 235.33 & 235.33 & 11.85 & 28.732 & 0.787 & 2.668 \\
& COLGEN 2 & 7.93 & 63.44 & - & 1.040 & 0.0448 & 0.082 \\
\hline
\end{tabular}
\caption{Column generation results for $d = 4$: number of iterations, number of columns generated, last pricing gap and running times (pricing model, restricted master model, hyperplane choice model).}
\label{tab:massive_colgen_4_1}
\end{table}

\begin{table}[h!]
\small
\centering
\begin{tabular}{c l r r r r r r}
$\mathcal{K}$ & method & \makecell{n. iterations} & \makecell{n. columns} & \makecell{last pricing gap} & \makecell[r]{PM time\\(seconds)} & \makecell[r]{RMM time\\(seconds)} & \makecell[r]{HCM time\\(seconds)} \\
\hline
\multirow{2}{0.2cm}{8} & COLGEN 1 & 210.88 & 210.88 & 2118.737 & 474.796 & 10.961 & 113.015 \\
& COLGEN 2 & 17.24 & 137.92 & - & 7.843 & 5.114 & 733.746 \\
\hline
\multirow{2}{0.2cm}{15} & COLGEN 1 & 182.15 & 182.15 & 1030.32 & 375.068 & 8.613 & 46.755 \\
& COLGEN 2 & 125.34 & 1002.72 & - & 54.187 & 19.545 & 528.832 \\
\hline
\multirow{2}{0.2cm}{16} & COLGEN 1 & 136.49 & 136.49 & 617.52 & 244.839 & 6.478 & 34.795 \\
& COLGEN 2 & 21.41 & 171.28 & - & 9.248 & 3.698 & 13.718 \\
\hline
\end{tabular}
\caption{Column generation results for $d = 8$: number of iterations, number of columns generated, last pricing gap and running times (pricing model, restricted master model, hyperplane choice model).}
\label{tab:massive_colgen_8_1}
\end{table}

\end{document}